\documentclass[dvips,11pt]{article}
\usepackage{amsmath,amssymb,amsfonts,amscd,rotating,pstricks,pst-node,bm,xspace}
\def\simarrow{\mathrel{\raise -0.5mm\hbox{$\sim$}}\hspace{-1.8mm}{\rightarrow} } 

\def\bsimarrow{\leftarrow\hspace{-0.7mm}\mathrel{\raise -0.5mm\hbox{$\backsim$}} }

\def\bt{\begin{tabular}}
\def\te{\end{tabular}}

\def\lettrine#1#2#3{\noindent\hangindent#1\hangafter-#2
\hskip-#1\smash{\hbox to #1{#3\hfill}}\ignorespaces}

\newcommand{\To}[1]{\mathop{\to}\limits_{#1}}

\def\BM{\begin{pmatrix}}
\def\EM{\end{pmatrix}}

\def\d=f{\buildrel\hbox{\scriptsize d\'{e}f}\over \Longleftrightarrow}

\def\square{\hfill\hbox{\vrule height .9ex width .8ex depth -.1ex}}
\def\rit{\text{\it I\hskip -2pt  R}}

\def\rl {\rit^{\hskip 1pt\ell}}

\def\Bd{{\text B}}

\def\Ed{{\text E}}

\def\Fs{{\cal F}}

\def\Hs{{\cal H}}

\def\Ms{{\cal M}}

\def\ung{\hbox{1\hskip -4.2pt \rm 1}}

\def\be{\begin{equation}}
\def\ee{\end{equation}}
\def\beqn{\begin{eqnarray}}
\def\eeqn{\end{eqnarray}}
\def\nobeqn{\begin{eqnarray*}}
\def\noeeqn{\end{eqnarray*}}
\def\ba{\left(\begin{array}}
\def\ea{\end{array} \right) }

\def\epr{\square\vskip 6pt}

\def\u{\underline}
\def\o{\overline}
\def\and{\; \mbox{and} \;}

\newcommand{\half}{\frac{1}{2}}

\def\hfl#1#2{\smash{\mathop{\hbox to 12mm{\rightarrowfill}}
\limits^{\scriptstyle #1}_{\scriptstyle #2}}}

\def\Be{\begin{enumerate}}
\def\Ee{\end{enumerate}}

\def\Bena{\begin{enumerate}
\def\labelenumi{\theenumi)}
\def\theenumi{\arabic{enumi}}
\def\labelenumii{\theenumii)}
\def\theenumii{\alph{enumii}}}

\def\Bean{\begin{enumerate}
\def\labelenumii{\theenumii)}
\def\theenumii{\arabic{enumii}}
\def\labelenumi{\theenumi)}
\def\theenumi{\alph{enumi}}}

\def\Bero{\begin{enumerate}
\def\labelenumii{\theenumii)}
\def\theenumii{\arabic{enumii}}
\def\labelenumi{(\theenumi)}
\def\theenumi{\roman{enumi}}}

\def\BeRo{\begin{enumerate}
\def\labelenumii{\theenumii)}
\def\theenumii{\arabic{enumii}}
\def\labelenumi{(\theenumi)}
\def\theenumi{\Roman{enumi}}}

\def\Bi{\vskip 11pt\begin{itemize}\itemsep=18pt}
\def\bi{\begin{itemize}}
\def\Ei{\end{itemize}\vskip 11pt}
\def\ei{\end{itemize}}

\def\Bd{\begin{description}}
\def\Ed{\end{description}}

\def\R{\right}
\def\L{\left}

\usepackage[latin1]{inputenc}

\def\sm{(semi)}
\def\ST{structures\xspace}
\def\smst{semistructure\xspace}
\def\SMST{semistructures\xspace}
\def\bsmst{bisemistructure\xspace}
\def\BSMST{bisemistructures\xspace}

\def\Times{{\mathrel{\u\times}}}
\def\gg{{\bm{g}}}

\def\prod{\mathop{\Pi}\limits}
\def\sum{\mathop{\Sigma}\limits}

\def\bbf{\boldmath\bf}
\def\o{\overline}

\def\Bi{\begin{itemize}}
\def\Ei{\end{itemize}}

\newcommand{\ZZ}{\mathbb{Z}\,}
\newcommand{\NN}{\mathbb{N}\,}
\newcommand{\PP}{\mathbb{P}\,}

\newcommand{\CC}{\mathbb{C}\,}

\def\res{\operatorname{res}}
\def\spec{\operatorname{spec}}

\def\End{\operatorname{End}}

\def\Hom{\operatorname{Hom}}
\def\corr{\operatorname{corr}}

\def\CH{\operatorname{CH}}

\def\GL{\operatorname{GL}}
\def\gl{\operatorname{gl}}

\def\bM{\begin{matrix}}
\def\eM{\end{matrix}}

\def\lr{left (resp. right) }
\def\rl{right (resp. left) }

\def\resp#1{(resp. #1)}
\def\rresp#1{\qquad \mbox{(resp.} \quad #1\ )}

\def\To{\begin{CD} @>>>\end{CD}}

\def\RL{_{R\times L}}
\def\Rl{_{R,L}}
\def\Lr{_{L,R}}

\def\wt{\widetilde}
\def\wh{\widehat}

\begin{document}

\setcounter{page}{0}
{\pagestyle{empty}
\null\vfill
\begin{center}
{\LARGE Introducing bisemistructures}
\vfill
{\sc C. Pierre\/}
\vskip 11pt

Institut de Mathématique pure et appliquée\\
Université de Louvain\\
Chemin du Cyclotron, 2\\
B-1348 Louvain-la-Neuve,  Belgium\\
pierre@math.ucl.ac.be

\vfill

\begin{abstract}
New fundamental mathematical structures are introduced by the triples (left semistructure, right semistructure, bisemistructure) associated with the classical mathematical structures and such that the bisemistructures, resulting from the reciprocal actions of left semistructures on right semistructures, are composed of bielements which are either diagonal bielements or cross bielements.

\end{abstract}
\vfill
\eject
\end{center}

\vfill\eject}
\setcounter{page}{1}
\def\thepage{\arabic{page}}
{\parindent=0pt 
\setcounter{section}{0}
\section*{Introduction}
\addcontentsline{toc}{section}{Introduction}

The existence of mathematical symmetric semiobjects (or semistructures), being for example semigroups, semirings or semimodules, is directly related to the generation of two symmetric left and right $n$-dimensional  affine semispaces $W_L$ and $W_R$, localized respectively in the upper and lower half spaces, from a polynomial ring $A=k[x_1,\dots,x_n]$ over a number field $k$ of characteristic $0$.

It is then shown in chapter 1 that {\bbf a mathematical structure can be split into two symmetric (disjoint) left and right semistructures referring (or localized in)  respectively to the upper and lower half spaces\/}.

But, a left and a right semistructure interact throughout the existence of a bisemistructure, which leads to consider {\bbf the triple (left semistructure, right semistructure, bisemistructure)\/}, associated with a given mathematical structure and verifying the following lemma:
\vskip 11pt

{\bf Lemma\/}: {\bfseries the right and left symmetric semistructures of a given structure operate on each other by means of their product giving rise to a bisemistructure in such a way that its bielements be either diagonal bielements or cross bielements.}
\vskip 11pt

Remark that a diagonal bielement can be defined with respect to a diagonal bilinear basis while an off-diagonal bielement can be developed according to an off-diagonal bilinear basis.
\vskip 11pt

{\bbf This concept of bi(semi)structure appears clearly in the endomorphism of a central simple algebra $B$ by its realization throughout the enveloping algebra $B^{\rm e}=B\otimes B^{\rm op}$, defined from the tensor product of this algebra $B$ by its opposite algebra $B^{\rm op}$.}
\vskip 11pt

It is also the case in the Lefschetz trace formula referring to the endomorphism of a complex of hypercohomology $R\Gamma  (X,L)$; where $L$ is a complex of sheaves on $X$, which can be calculated according to classes of cohomology on $X\times X$ \cite{G-I}.

In this perspective, it is perhaps not innocent to remark that {\bbf a matrix is basically a bilinear object\/} and that the representation space of an algebraic group of matrices $\GL(n,\rit)$ is a vectorial space of dimension $n^2$ \cite{B-T}, \cite{D-G}, \cite{Bou}.
\vskip 11pt

{\bbf Fundamental semistructures and bisemistructures\/}, then {\bf introduced\/} and defined {\bf in chapter 2\/}, are:
\Bi
\item left and right semigroups and bilinear semigroups.
\item left and right monoids and bimonoids.
\item left and right semirings and bisemirings.
\item left and right semifields and bisemifields.
\item left and right semimodules and bisemimodules.
\item inner product bisemispaces, among which bilinear Hilbert semispaces.
\item left and right semialgebras and bisemialgebras.
	   \Ei
	   
	   {\bbf The fundamental triple is $(G_L,G_R,G\RL)$\/} where:
	   \Bi
	   \item {\bbf $G_L$ \resp{$G_R$} is a \lr semigroup under the addition\/} of its \lr elements $g_{L_i}$ \resp{$g_{R_i}$} restricted to (or referring to) the upper (resp. lower) half space.
	   
	   \item {\bbf $G\RL$ is a bisemigroup, or a bilinear semigroup, such that its bielements\linebreak $(g_{R_i}\times g_{L_i})$ be submitted to the cross binary operation $\Times$ according to:\/}
	   \begin{align*}
	  G\RL \Times G\RL
	  &\To G\RL\\
	    (g_{R_i}\times g_{L_i}) \Times (g_{R_j}\times g_{L_j})
	  &\To (g_{R_i}+ g_{R_j}) {\times} (g_{L_i}+ g_{L_j})\end{align*}
	  {\bbf leading to cross products\/} 
	  $(g_{R_i}\times g_{L_j}) $ and $(g_{R_j}\times g_{L_i})$.
	  \Ei
	  \vskip 11pt
	  
	  Thus,  the cross binary operation $\Times$, responsible for the bilinear character of $G\RL$, generalizes the concept of (semi)group $G\RL$ by sending pairs of right and left elements either in diagonal bielements $(g_{R_i}\times g_{L_i}) $ and $(g_{R_j}\times g_{L_j})$ or in cross bielements $(g_{R_i}\times g_{L_j}) $ and $(g_{R_j}\times g_{L_i})$.
	  
	  The bilinear semigroup $G\RL=G_R\times G_L$ is thus a generalization of the direct product $G_R\times G_L$ of the (semi)groups $G_R$ and $G_L$ since, in this case, only diagonal bielements are generated.
	  \vskip 11pt
	  
	 Another important triple is $(M_L,M_R,M \RL)$ where:
	 \Bi
	 \item $M_L$ \resp{$M_R$} is a \lr $R_L$-semimodule $M_L$ \resp{$R_R$-semimodule $M_R$} over the \lr semiring $R_L$ \resp{$R_R$}.
	 
	 \item $M\RL=M_R\times M_L\simeq M_R\otimes_{R_R\times R_L} M_L$ is a $R_R\times R_L$-bisemimodule over the bisemiring $R\RL=R_R\times R_L$, i.e. an additive abelian bisemigroup.
	 \Ei
	 \vskip 11pt
	 
	{\bbf The $R\RL$-bisemimodule $M_R\otimes_{R_R\times R_L} M_L$ splits naturally into:
	\[
	M_R\otimes_{R\RL} M_L=(M_R\otimes _D M_L)+(M_R\otimes_{OD} M_L)\]}
	where:
	\Bi
	\item {\bbf $(M_R\otimes_D M_L)$ is a diagonal $R\RL$-bisemimodule\/} of which bielements are expressed in a diagonal bilinear basis $\{e_\alpha \otimes f_\alpha \}_\alpha $, $e_\alpha \in M_R$, $f_\alpha \in M_L$.
	
	\item {\bbf $(M_R\otimes_{OD}M_L)$ is an off-diagonal $R\RL$-bisemimodule\/} over $R\RL$ of which bielements are expressed in an off-diagonal bilinear basis $\{e_\alpha \otimes f_\beta  \}_{\alpha\neq\beta } $.
	\Ei
	\vskip 11pt
	
This leads us to introduce {\bbf mixed product $R_R\times R_L$-bisemispaces $M_R\otimes M_L$\/} endowed with bilinear functions $[\centerdot,\centerdot]$ from $M_R\times M_L$ into $R_R\times R_L$: they can be extended, diagonal or off-diagonal mixed product bisemispaces $M_R\otimes M_L$, $M_R\otimes_D M_L$ or
$M_R\otimes_{OD} M_L$ according as the considered linear basis is complete, diagonal or off-diagonal.

These mixed product bisemispaces become {\bbf external \lr product bisemi\-spaces\/} by the projection of the $R_R$-semispace $M_R$
\resp{$R_L$-semispace $M_L$} onto the
$R_L$-semispace $M_L$
\resp{$R_R$-semispace $M_R$}.  And, if every covariant \resp{contravariant} element of the projected $M_R$ \resp{$M_L$} onto $M_L$ \resp{$M_R$} is mapped onto the corresponding contravariant \resp{covariant} element, then these external \lr product bisemispaces {\bbf are transformed into corresponding \lr inner product bisemispaces\/}. It is then shown that the \lr diagonal inner product bisemispace is a \lr separable bilinear Hilbert semispace being in one-to-one correspondence with the classical separable linear Hilbert \sm space.
\vskip 11pt

The last proposed fundamental triple is $(A_L,A_R,A_R\otimes A_L)$ where:
\Bi
\item $A_L$ \resp{$A_R$} is a \lr $R_L$-semialgebra
\resp{$R_R$-semialgebra} over the \lr (division) semiring $R_L$ \resp{$R_R$}.

\item {\bbf $A_R\otimes A_L$ is a $R_R\times R_L$-bisemialgebra\/} which is a bisemiring $A\RL$ such that the pair $(A\RL,+)$ is an unitary $R\RL$-bisemimodule equipped with the bilinear homomorphism
\[ \mu \Rl: \quad A\RL \Times A\RL\To A\RL\]
where $\Times$ is the cross binary operation allowing to generate off-diagonal bielements.
\Ei

This enables to introduce {\bbf a bisemialgebra of Hopf\/} given by the triple $((A_{R(P)}\otimes A_L),\linebreak 
(A_{L(P)}\otimes A_R),S_b)$ where:
\Bean
\item $(A_{R(P)}\otimes A_L)$ \resp{$(A_{L(P)}\otimes A_R)$} is a \lr bisemialgebra in such a way that $A_{R(P)}$ \resp{$A_{L(P)}$} is the \rl semialgebra $A_R$ \resp{$A_L$} projected onto 
$A_{L}$ \resp{$A_R$}.

\item $S_b:A_{R(P)}\otimes A_L \to A_{L(P)}\otimes A_R$ is the antipode mapping bijectively the left bisemialgebra $A_{R(P)}\otimes A_L$ onto the right bisemialgebra $A_{L(P)}\otimes A_R$.
\Ee

Working with bisemialgebras thus allows to define very clearly the antipode $S_b$ and its inverse $S_b^{-1}$.
\vskip 11pt

Chapter 3 is devoted to the {\bbf study of\/} concrete bisemistructures, under the circumstances {\bbf bilinear algebraic semigroups and bisemischemes\/}.
\vskip 11pt

{\bf In conclusion\/}, I would like to emphasize that {\bbf \BSMST are richer than the corresponding classical \ST by the existence of cross bisemisubstructures\/}.
\vskip 11pt

On the other hand, {\bbf bielements arise naturally\/}.  For example, consider the conjugation action of a subgroup $H$ on a group $G$ given by $(h,x)\to h\mathrel x h^{-1}$, $\forall\ h\in H$, $x\in G$.  In the case of \BSMST, we would have to consider the conjugation biaction of a subbisemigroup $H_R\times H_L$ on the bisemigroup $G_R\times G_L$ given by $(h_R\times h_L, x_R\times x_L)\to h_L\ x_L\ x_R\ h_R$, where $h_R=h_L^{-1}$, $x_L\in g_L$, $x_R\in g_R$, $h_L\in H_L$, $h_R\in H_R$.
\vskip 11pt

In this philosophy of \BSMST, the {\bbf left \smst acts on the corresponding right \smst\/} and vice versa {\bbf giving a better visibility to their endomorphisms\/} and to some \BSMST involving internal structural endomorphisms, as for example the antipode map in the (bisemi)algebra of Hopf.
\vskip 11pt

Finally, in this context, the metric, the scalar product, the norm, the invariant bilinear forms of Lie groups \cite{Sch3},\dots\  get a new ontological meaning: this is tied up to the fact that, {\bbf in a given (bisemi)structure, only the left semistructure is generally observable in such a way that the (tensor) product of this left semistructure by the opposite right (dual) hidden semistructure constitutes the only mathematical manageable (bisemi)structure corresponding to the objectivable reality\/} \cite{Pie2}.

\section{Origin of \SMST and \BSMST}

\subsection{Origin of semiobjects}

{\bbf The existence of mathematical symmetric semiobjects (or \SMST) results essentially from classical algebraic geometry\/}.
\vskip 11pt

Let $A=k[x_1,x_2,\dots,x_n]$ be a polynomial ring at $n$ indeterminates over a global number field $k$ of characteristic $0$.  A solution of the system of polynomial equations $P_\mu (x_1,\dots,x_n)=0$ of the ideal $I_L$ of $k[x_1,\dots,x_n]$ is an $n$-tuple $(a_1,\dots,a_n)$.

But the existence of the ideal $I_R$ of $k[x_1,\dots,x_n]$, composed of polynomial equations\linebreak $P_\mu (-x_1,\dots,-x_n)=0$, provides solutions which are negative $n$-tuples $(-a_1, \dots,-a_n)$.  Thus, the set of zeros of $I_L\cup I_R=\{P_\mu (x_1,x_2,\dots,x_n)\cup P_\mu (-x_1,\dots,-x_n)\mid P_\mu (W)=0\}$ is an affine $n$-dimensional space $W$ which can be split into two symmetric semispaces $W_L$ and $W_R$ in such a way that:
\Bi
\item $W=W_L\cup W_R$; $W_L\cap W_R=\emptyset$.

\item $W_L$ is an affine $n$-dimensional semispace restricted to the upper half space.

\item $W_R$ is the affine $n$-dimensional semispace, symmetric to $W_L$, restricted to the lower half space and disjoint of $W_L$.
\Ei

So, the consideration of a general polynomial ring $k[x_1,\dots,x_n]$ allows to generate the union $W_R\cup W_L$ of two symmetric right and left $n$-dimensional affine semispaces $W_R$ and $W_L$.
\vskip 11pt

\subsection{On the necessity of considering bisemiobjects}

On the other hand, in the case of number fields, the ``total'' (closed) algebraic extension of the number field $k$ will be the union $\wt  F=\wt  F_R\cup \wt  F_L$ of two symmetric (closed) algebraic extensions $\wt  F_R$ and $\wt  F_L$ referring, for example, respectively to the sets of negative and positive roots of the polynomial ring $k[x]$.
\vskip 11pt

So, taking into account that:
\Be
\item a left $T_n(\wt  F_L)$-semimodule $T^{(n)}(\wt F_L)$ \resp{right $T^t_n(\wt  F_R)$-semimodule $T^{(n)}(\wt F_R)$} is generated under the \lr action of the upper (resp. lower) triangular group of matrices 
$T_n(\wt F_L)$
\resp{$T_n^t(\wt F_R)$} according to \cite{Pie1}:
\begin{align*}
T_n(\centerdot): \quad \wt  F_L &\To T^{(n)}(\wt  F_L)\;, && \text{also noted $V_L$}\\
\rresp{T_n^t(\centerdot): \quad \wt  F_R &\To T^{(n)}(\wt  F_R)\;, && \text{also noted $V_R$}}.\end{align*}

\item {\bbf the (bilinear) algebraic (semi)group of matrices\/}:
\[ \GL_n(\wt  F_R\times \wt  F_L)=T^t_n(\wt  F_R)\times T_n(\wt  F_L)\;, \]
written in function of the Gauss bilinear decomposition, {\bbf generates the affine bisemispace\/}
$T^{(n)}(\wt  F_R)\otimes_{\wt  F_R\times \wt  F_L} T^{(n)}(\wt  F_L)
\equiv V_R \otimes_{\wt  F_R\times \wt  F_L} V_L$
(as it will be developed subsequently),
\Ee

it appears clearly that a bisemiobject (or a \bsmst), resulting from the (tensor) product of a left semiobject (or \smst) by its symmetric right equivalent, must also be envisaged: this leads us to consider triples: $(V_R,V_L,V_R
\otimes_{\wt  F_R\times \wt  F_L} V_L)$ and, more generally, {\bbf triples of \ST\/} (right \smst, left \smst, \bsmst) because right and left symmetric \SMST ``interact'' by means of \BSMST, as stated in the following lemma:
\vskip 11pt

{\bf Lemma\/}: {\em the right and left symmetric semistructures of a given structure operate always on each other by means of their product giving rise to a bisemistructure in such a way that its bielements be either diagonal bielements or cross bielements.}
\vskip 11pt

\paragraph{Proof}: Let us only say that a \bsmst is basically a bilinear semigroup $G\RL$ (as it will be developed subsequently), such that:
\Bean
\item the right and left symmetric \SMST are semigroups $G_R$ and $G_L$ under the addition of their respective elements.

\item the bielements of this bilinear semigroup split into diagonal bielements or cross bielements according to the cross binary operation $\Times$ of $G\RL$ (see definition 2.1.2).\epr
\Ee
\vskip 11pt

Let us remark that the consideration of triple of \ST reminds {\bbf the construction of pure motives\/} on algebraic varieties embodied in their groups of classes of cycles \cite{Gro}.  Indeed, a pure (Chow) motive is essentially a pair $(X,p)$ where $X$ is a $n$-dimensional smooth projective variety and $p\in\corr^0(X,X)=\CH^{2n}(X\times X)\simeq\dfrac{Z^i(X\times X)}{Z^i_{{\rm rat}}(X\times X)}$ is the projector of the Chow ring\linebreak $\CH^{2n}(X\times X)$ of the group $Z^i(X\times X)$ of algebraic (bi)cycles on $(X\times X)$ by the subgroup $Z^i_{{\rm rat}}(X\times X)$ rationally equivalent to $0$ \cite{Mur}.
\vskip 11pt

{\bbf Most of bisemiobjects, which are (bilinear) products of semiobjects, proceed generally from the bilinear structure of matrices\/}.  In this respect, every endomorphism of a central simple algebra $B$ can only be handled throughout his enveloping algebra $B^e=B\otimes_{\wt  F}B^{{\rm op}}$, where $B^{\rm op}$ denotes the opposite algebra of $B$; indeed, in that case, we have that:
\begin{equation}\label{eq:*}
B^e\approx \End_{\wt  F}(B)
 \approx \End_{\wt  F}(\wt F^n)
 \approx \Ms_n(\wt  F)\tag{$*$}\end{equation}
 where $\Ms_n(\wt  F)$, being the set of $(n\times n)$ matrices over the field $\wt  F$, is a central simple algebra over  $\wt  F$ \cite{F-D}.
 \vskip 11pt
 
 If an object is then envisaged as a bisemiobject, it becomes clear from the preceding reflection on endomorphism that {\bbf every endomorphism of a (semi)object\/}, describing the variations of the internal algebraic (and/or geometric) structure of this one, {\bbf can only be realized throughout the (tensor) product of this semiobject with its (``opposite'') symmetric equivalent\/}.
 \vskip 11pt
 
 Indeed, the isomorphisms \eqref{eq:*}, transposed in the case of a central simple semialgebra $A$ over $\wt  F_L$ (which will be defined in a next section) associated with  left semiobjects, become:
 \[ A^e=A\otimes_{\wt  F_R\times \wt  F_L} A^{\rm op}
 \simeq \End_{\wt  F_R\times \wt  F_L}(A\times A^{\rm op})
 \simeq \End_{\wt  F_R\times \wt  F_L}((\wt  F_R\times \wt  F_L)^n)
 \simeq \GL_n(\wt  F_R\times \wt  F_L)\]
 with the evident notations.
 \vskip 11pt
 
 \subsection{Definitions: algebraic bisubsets and bipoints}
 
 Assume that the two symmetric \SMST are the right and left $n$-dimensional affine semispaces $V_R$ and $V_L$.
 
 \Bi
 \item {\bbf A simple \lr  subset\/} of $V_L$ \resp{$V_R$} {\bf is}:
 \Be
 \item {\bbf a one-dimensional irreducible closed algebraic subset\/} $\wt F^+_{v^1_j}$ \resp{$\wt F^+_{\o v^1_j}$} (see notations of \cite{Pie1}) characterized by the Galois extension degree:
 \[ [\wt F^+_{v^1_j}:k]=
[\wt F^+_{\o v^1_j}:k]=N\]
which is an algebraic dimension.

\item a simple (irreducible) \lr $k$-semimodule, interpreted in \cite{Pie2} as a \lr algebraic quantum.
\Ee

\item {\bbf A  simple bisubset\/} of $V_R\otimes_{\wt F_R\times \wt F_L} V_L$ {\bbf is an irreducible closed algebraic bisubset\/} $(\wt F^+_{\o v^1_j}\times_D \wt F^+_{v^1_j})$ characterized by the Galois extension bidegree (which is an algebraic dimension):
\[ [\wt F^+_{\o v^1_j}:k] \times_D [\wt F^+_{v^1_j}:k]=N\]
referring to a ``diagonal'' product ($\times _D$) where the off-diagonal components are not taken into account.

\item An element of an irreducible closed algebraic \lr subset
$\wt F^+_{v^1_j}\in \wt F_L$
\resp{$\wt F^+_{\o v^1_j}\in \wt F_R$} is {\bbf an algebraic \lr point\/} 
$(0,\dots,a_k,\dots,0)$
\resp{$(0,\dots,-a_k,\dots,\linebreak 0)$} or a component of a \lr point of $V_L$ \resp{$V_R$}.

\item A diagonal bielement of an irreducible closed algebraic bisubset $(\wt F^+_{\o v^1_j}\times_D \wt F^+_{v^1_j})$ is {\bf an algebraic diagonal bipoint\/}
$(0,\dots,(-a_k)\times_D (a_k),\dots,0)$ or a component of a diagonal bipoint of $V_R\otimes V_L$.

\item An off-diagonal bielement of an off-diagonal closed algebraic bisubset
$(\wt F^+_{\o v^1_j}\times_{OD} \wt F^+_{v^1_j})$ is {\bf an algebraic off-diagonal bipoint\/}
$(0,\dots,(-b_k)\times_{OD} (c_k),\dots,0)$
characterized by a nonzero bilinear ``off-diagonal'' component (from which the notation ``$\times_{OD}$'' refers to an off-diagonal product).

Remark that
\[ (\wt F^+_{\o v^1_j}\times \wt F^+_{v^1_j})
= (\wt F^+_{\o v^1_j}\times_D \wt F^+_{v^1_j})
+ (\wt F^+_{\o v^1_j}\times_{OD} \wt F^+_{v^1_j})\]
as it will be shown in the following.

\item {\bbf A diagonal algebraic bipoint\/} of the affine bisemispace $(V_R\otimes_{\wt F_R\times \wt F_L} V_L)$ is given by the $n$-bituple $(-a_1\times_D a_1,\dots,
-a_k\times_D a_k,\dots,
-a_n\times_D a_n)$ while a {\bbf a general algebraic bipoint\/} of $(V_R \otimes_{\wt F_R\times \wt F_L} V_L)$ is described by the $n^2$-bituple 
\[\bM
(-a_1\times_D a_1,\dots,-a_1\times_{OD} a_n,\\[-6pt]
\ -a_2\times_{OD} a_1,\dots,-a_2\times_{OD} a_n,\\[-6pt]
\vdots\\[-6pt]
\ -a_n\times_{OD} a_1,\dots,-a_n\times_{D} a_n)\;.\eM\]
\Ei\vskip 11pt

\subsection{Definitions: geometric bipoints}

\Bi
\item {\bbf A general geometric bipoint\/} of $(V_R\otimes_{\wt F_R\times \wt F_L} V_L)$ {\bf is a general algebraic bipoint given by the $n^2$-bituple\/}
$(-a_1\times_D a_1,\dots,-a_i\times_{OD}a_j,\dots,-a_n\times_D a_n)$, $1\le i,j\le n$,
{\bbf together with the additional geometric $n^2$-bituple\/} $(g_{11},\dots,g_{ij},\dots,g_{nn})$ where the $g_{ij}$ are the components at this bipoint of the metric tensor $\gg$ of type $(0,2)$ with respect to a bilinear basis $\{(e_i\otimes e_j)\}_{i,j}$ in such a way that the
$g_{ij}=g(e_i,e_j)$ are the scalar products of the basis vectors.

\item On the other hand, {\bbf if the right semispace $V_R$ is projected onto the left semispace $V_L$\/} generating the bisemispace $(V_{R_P}\otimes_{\wt F_R\times \wt F_L}V_L)$, {\bbf then a general geometric bipoint of 
$(V_{R_P}\otimes_{\wt F_R\times \wt F_L}V_L)$ will be given by\/}:
\Bena
\item {\bbf the algebraic $n^2$-bituple\/}:
$(a_1\times_D a_1,\dots,a_i\times_{OD}a_j,\dots,a_n\times_D a_n)$ where:
\[ \PP_{R\to L}(i): \quad -a_i\To +a_i\;, \qquad \forall\ 1\le i\le n\;, \]
is the map projecting the $i$-th right coordinate $-a_i$ of the general algebraic bipoint of $(V_{R}\otimes_{\wt F_R\times \wt F_L}V_L)$ into the $i$-th right coordinate $+a_i$ of the corresponding algebraic bipoint on
$(V_{R_P}\otimes_{\wt F_R\times \wt F_L}V_L)$.

\item {\bbf the additional geometric $n^2$-bituple\/} $(g^1_1,\dots,g^j_i,\dots g^n_n)$ where the $g^j_i$ are the components at this bipoint of the metric tensor of type $(1,1)$ with respect to the bilinear basis $\{(e_i\otimes e^j)\}_{i,j}$ in such a way that 
$g_{i}^j=g(e_i,e^j)$.
\Ee

{\bbf This geometric bipoint\/} of
$(V_{R_P}\otimes_{\wt F_R\times \wt F_L}V_L)$
{\bbf can be identified with a classical geometric point of $V_L$\/} which is characterized by:
\Bena
\item the algebraic $n$-tuple $(a_1,\dots,a_j,\dots,a_n)$.
\item the additional geometric $n^2$-bituple $(g_{11},\dots,g_{ij},\dots,g_{nn})$ of which components are those of the symmetric metric tensor of type $(0,2)$.
\Ee

\item Similarly, if the left semispace 
{  $V_L$ is projected onto  $V_R$\/} generating the bisemispace $(V_{R}\otimes_{\wt F_R\times \wt F_{L_P}}V_{L_P})$, { then a general geometric bipoint of 
$(V_{R}\otimes_{\wt F_R\times \wt F_{L_P}}V_{L_P})$ will be given by\/}:
\Bena
\item {\bbf the algebraic $n^2$-bituple\/}:
$((-a_1)\times_D (-a_1),\dots,(-a_i)\times_{OD}(-a_j),\dots, (-a_n)\times_D (-a_n))$.
\item the additional geometric $n^2$-bituple $(g^1_{1},\dots,g^i_{j},\dots,g^n_{n})$ where the $g^i_j$ are the components of the  metric tensor of type $(1,1)$ at this bipoint.
\Ee
\Ei

\section{Fundamental \SMST and \BSMST}

Right and left \SMST and the corresponding \BSMST will now be recalled or introduced.  The notation $L,R$ means ``left'' or ``right''.
\vskip 11pt

\subsection[The bilinear semigroup $G\RL$]{\bbf The bilinear semigroup $G\RL$}

\subsubsection{Definition}

{\bbf A \lr semigroup $G_L$ \resp{$G_R$}\/}, noted in condensed notation $G\Lr$, is a nonempty set of \lr elements $g_{L_i}$ \resp{$g_{R_i}$}, localized in (or referring to) the upper \resp{lower} half space, together with a binary operation (addition or multiplication) on $G\Lr$, i.e. a function $G\Lr\times G\Lr \to G\Lr$ (mostly, additive case) or $G\Lr\times G_{L,L} \to G\Lr$ (mostly, multiplicative case) which is associative:
\begin{align*}
(g_{L_i}\times g_{L_j})\times g_{L_k}
&=g_{L_i}\times (g_{L_j}\times g_{L_k})\;, \qquad g_{L_i},g_{L_j},g_{L_k}\in G_L\\
\rresp{(g_{R_i}\times g_{R_j})\times g_{R_k}
&=g_{R_i}\times (g_{R_j}\times g_{R_k})\;, \qquad g_{R_i},g_{R_j},g_{R_k}\in G_R}.\end{align*}
\vskip 11pt

\subsubsection{Definition}

{\bbf The bilinear semigroup $G\RL$\/}, also called bisemigroup, {\bbf referring to the semigroups $G_R$ and $G_L$, is defined by the bilinear function\/}
\begin{align*}
 G_R\times G_L&\To G\RL\\
(g_{R_i},g_{L_i}) &\To (g_{R_i}\times g_{L_i})\end{align*}
\Bean
\item[$\alpha$)] sending pairs $(g_{R_i},g_{L_i})$ of symmetric elements, localized in (or referring to) the lower and upper half spaces, to their products
$(g_{R_i}\times g_{L_i})$.

\item[$\beta$)] {\bbf submitted to the cross binary operation ``$\Times$''}: 
\begin{multline*}
G\RL\Times G\RL\to G\RL\ \text{defined by:}\quad (g_{R_i}\times g_{L_i})\Times (g_{R_j}\times g_{L_j})\\
\To (g_{R_i}+g_{R_j})
\times (g_{L_i}+g_{L_j})\end{multline*}
{\bbf in such a way that cross products\/}
$(g_{R_i}\times g_{L_j})$ and
$(g_{R_j}\times g_{L_i})$, resulting from the bilinearity of $G\RL$, {\bf could be generated\/}.

Indeed, the development of
$(g_{R_i}+g_{R_j})
\times (g_{L_i}+g_{L_j})$ gives
\[
(g_{R_i}+g_{R_j})\times (g_{L_i}+g_{L_j})
=(g_{R_i}\times g_{L_i})+ (g_{R_j}\times g_{L_j})
+(g_{R_i}\times g_{L_j})+ (g_{R_j}\times g_{L_i})\]
where $(g_{R_i}\times g_{L_i})$
and $(g_{R_j}\times g_{L_j})$ are ``diagonal'' bielements.
\Ee

But, as the elements $g_{R_i}$, $g_{L_i}$, $g_{R_j}$ and $g_{L_j}$ are simple (i.e. non composite), they can only appear once in the cross binary operation
$(g_{R_i}\times g_{L_i})\Times (g_{R_j}\times g_{L_j})$. 

 So, {\bbf a left or right element in 
$(g_{R_i}+g_{R_j})\times (g_{L_i}+g_{L_j})$ appears\/}:
\Bean
\item {\bf either in\/} the diagonal products
$(g_{R_i}\times g_{L_i})$ and $ (g_{R_j}\times g_{L_j})$ of right and left elements, i.e. in {\bf ``diagonal'' bielements\/}.

\item {\bbf or in\/} the (off-diagonal or) cross products 
$(g_{R_i}\times g_{L_j})$ and $ (g_{R_j}\times g_{L_i})$ of right and left elements, i.e. in {\bbf ``cross'' bielements\/}.
\Ee
\vskip 11pt

This leads us to formulate the following propositions:
\vskip 11pt

\subsubsection{Proposition}

{\em The multiplication between bielements is the cross binary operation ``$\Times$''.}
\vskip 11pt

\paragraph{Proof.} \quad This results from definition 2.1.2.\epr
\vskip 11pt

\subsubsection{Proposition}

{\em The cross binary operation ``$\Times$'', responsible for the bilinear character of the semigroup $G\RL$, generalizes the concept of semigroup $G\RL$ by sending pairs of right and left elements either in diagonal bielements of in cross bielements.}
\vskip 11pt

\paragraph{Proof.} \quad The statement of this proposition results from:
\Bean
\item The fact that the cross binary operation allows to generate cross products;
\item the above mentioned properties of bielements of the semigroup $G\RL$.\epr
\Ee
\vskip 11pt

\subsubsection{Proposition}

{\em A bilinear semigroup $G\RL$ is abelian if its cross binary operation is commutative.}
\vskip 11pt

\paragraph{Proof.} \quad The cross binary operation is commutative if the bielements commute according to:
\[
(g_{R_i}\times g_{L_i})\Times (g_{R_j}\times g_{L_j})
=(g_{R_j}\times g_{L_j})\Times (g_{R_i}\times g_{L_i})\;.\]
Developing the left hand side, we get:
\begin{align*}
(g_{R_i}\times g_{L_i})\Times (g_{R_j}\times g_{L_j})
&=(g_{R_i}+ g_{R_j})\times (g_{L_i}+g_{L_j})\\
&=(g_{R_i}\times g_{L_i})+ (g_{R_j}\times g_{L_j})
+(g_{R_j}\times g_{L_i})+ (g_{R_i}\times g_{L_j})
\end{align*}
by taking into account the definition 2.1.2.

On the other hand, the development of the right hand side gives similarly:
\[
(g_{R_j}\times g_{L_j})\Times (g_{R_i}\times g_{L_i})
 =(g_{R_j}\times g_{L_j})+ (g_{R_i}\times g_{L_i})
+(g_{R_j}\times g_{L_i})+ (g_{R_i}+g_{L_j})\;.
\]
It then appears that a bilinear semigroup $G\RL$ is generally abelian because the abelian character of the cross binary operation, sending diagonal bielements either into themselves or into cross bielements, is reduced to the abelian character of the addition.\epr
\vskip 11pt

\subsubsection{Algebraic general and cross bipoints}

Let $(-a_{i_1},\dots,-a_{i_k},\dots,-a_{i_n})$ and
  $(-a_{j_1},\dots,-a_{j_k},\dots,-a_{j_n})$ be the coordinates of two points $a_{i_R}$ and $a_{j_R}$ of the $n$-dimensional right affine semivariety $V_R$ and let
  $(+a_{i_1},\dots,\linebreak +a_{i_k},\dots,+a_{i_n})$ and
  $(+a_{j_1},\dots,+a_{j_k},\dots,+a_{j_n})$ the coordinates of two symmetric points $a_{i_L}$ and $a_{j_L}$ in the left affine semivariety $V_L$.
  
  Then, the cross binary operation over the affine bisemivariety $(V_R\otimes_{F^0} V_L)$, which is fundamentally a bilinear semigroup, is given on the bipoints
  $(a_{i_R}\times a_{i_L})$ and
  $(a_{j_R}\times a_{j_L})$ by:
  \begin{align*}
  (a_{i_R}\times a_{i_L})\Times
  (a_{j_R}\times a_{j_L})
&= (a_{i_R}+ a_{j_R})\times(a_{i_L}+ a_{j_L})\\
&= (a_{i_R}\times a_{i_L})+(a_{j_R}\times a_{j_L})
+[(a_{i_R}\times a_{j_L})+(a_{j_R}\times a_{i_L})]
\end{align*}
where:
\Bi
\item $(a_{i_R}\times a_{i_L})$ and
  $(a_{j_R}\times a_{j_L})$ are {\bbf general algebraic bipoints\/} in such a way that the coordinates of
$(a_{i_R}\times a_{i_L})$ are the $n^2$-bituple $(-a_{i_1}\times_D a_{i_1},\dots,-a_{i_1}\times_{OD} a_{i_n},
\dots,\linebreak
-a_{i_n}\times_{OD} a_{i_1},\dots,-a_{i_n}\times_D a_{i_n})$;

\item $(a_{i_R}\times a_{j_L})$ and
  $(a_{j_R}\times a_{i_L})$ are {\bbf algebraic cross bipoints\/} in such a way that the cross coordinates of
  $(a_{i_R}\times a_{j_L})$ are the $n^2$-bituple
  $(-a_{i_1}\times_D a_{j_1},\dots,-a_{i_1}\times_{OD} a_{j_n},\dots,
-a_{i_n}\times_{OD} a_{j_1},\dots,-a_{i_n}\times_D a_{j_n})$.
\Ei
\vskip 11pt

\subsection{Bimonoids, bisemirings and bisemifields}

\subsubsection{Definition}

{\bbf A \lr monoid\/} is a \lr semigroup $G_L$ \resp{$G_R$} which contains an identity element $e_L$ \resp{$e_R$} such that:
\[ e_L\centerdot g_L=g_L\quad 
\rresp{ g_R \centerdot e_R=g_R}, \qquad \forall\ g_L\in G_L\ ,\;
g_R\in G_R\;.\]
Generally, $e_L=e_R$.
\vskip 11pt

A {\bf bimonoid\/} is a bisemigroup in which the identity element $e\RL=e_R\times e_L$ verifies: $e\RL=e_R=e_L$.
\vskip 11pt

\subsubsection{Definition}

{\bbf A \lr semiring\/} is a nonempty set $R_L$ \resp{$R_R$}, also written $R\Lr$, together with two binary operations (addition and multiplication) such that:
\Bena
\item $(R\Lr,+)$ is an abelian \lr semigroup (or monoid).

\item $(g_{L_i,R_i}\times g_{L_j,R_j})\times g_{L_k,R_k}=g_{L_i,R_i}\times (g_{L_j,R_j}\times g_{L_k,R_k})$: associative multiplication.

\item $\begin{aligned}[t]
& g_{L_i,R_i}\times (g_{L_j,R_j}+ g_{L_k,R_k})
=(g_{L_i,R_i}\times g_{L_j,R_j})+ (g_{L_i,R_i}\times g_{L_k,R_k})\\
\text{and\ }
& (g_{L_i,R_i}+ g_{L_j,R_j})\times g_{L_k,R_k}
=(g_{L_i,R_i}\times g_{L_k,R_k})+ (g_{L_j,R_j}\times g_{L_k,R_k})\ :\end{aligned}$

 left and right distribution.
\Ee

A {\bf   bisemiring\/} $R\RL$ is an abelian bisemigroup (or bimonoid) submitted to the two binary operations (addition and multiplication) verifying:
\Bena
\item associative multiplication (or cross binary product):
\[
[(g_{R_i}\times g_{L_i})
\centerdot (g_{R_j}\times g_{L_j})]
\centerdot (g_{R_k}\times g_{L_k})
=(g_{R_i}\times g_{L_i})
\centerdot [(g_{R_j}\times g_{L_j})
\centerdot (g_{R_k}\times g_{L_k})]\]
where the multiplication ``$\centerdot$'' is the cross binary operation $\Times$ according to definition 2.1.2.

So, we get:
\[
[(g_{R_i}\times g_{L_i})
\Times (g_{R_j}\times g_{L_j})]
\Times (g_{R_k}\times g_{L_k})
=(g_{R_i}\times g_{L_i})
\Times [(g_{R_j}\times g_{L_j})
\Times (g_{R_k}\times g_{L_k})]\]
which gives:
\[
[(g_{R_i}+ g_{R_j})+ g_{R_k}]
\times[(g_{L_i}+ g_{L_j})+ g_{L_k}]
=[g_{R_i}+ (g_{R_j}+ g_{R_k})]
\times[g_{L_i}+ (g_{L_j}+ g_{L_k})]\;.\]

\item left and right distribution with respect to cross binary product
\begin{multline*}
(g_{R_i}\times g_{L_i})
\centerdot [(g_{R_j}\times g_{L_j})
+ (g_{R_k}\times g_{L_k})]
\\ =[(g_{R_i}\times g_{L_i})
\centerdot (g_{R_j}\times g_{L_j})]
+ [(g_{R_i}\times g_{L_i})
\centerdot (g_{R_k}\times g_{L_k})]
\end{multline*}
which gives:
\begin{multline*}
(g_{R_i}\times g_{L_i})
\Times [(g_{R_j}\times g_{L_j})
+ (g_{R_k}\times g_{L_k})]\\
=[(g_{R_i}\times g_{L_i})
\Times (g_{R_j}\times g_{L_j})]
+ [(g_{R_i}\times g_{L_i})
\Times (g_{R_k}\times g_{L_k})]
\end{multline*}
if ``$\centerdot$'' is identified to ``$\Times$''.
\Ee
\vskip 11pt

\subsubsection{Definitions}

\Bean
\item A {\bf \lr integral domain\/} is a commutative \lr semiring $R\Lr$ with identity $\ung_{R\Lr}$ and no zero divisors.

An {\bf integral bidomain\/} is a commutative bisemiring $R\RL$ with identity $\ung_{R\RL}$ and no zero bidivisors (i.e. products, right by left, of divisors).

\item A {\bf \lr division semiring\/} $R\Lr$ is a \lr integral domain if every element of $R\Lr$ is a unit (i.e. invertible).

A {\bf division bisemiring\/} is an integral bidomain if every bielement of $R\RL$ is a biunit (i.e. product of a right unit by a left unit).

\item A {\bf \lr semifield\/} $\wt F_L$ \resp{$\wt F_R$} is a commutative \lr division semiring.

A {\bf bisemifield\/} $\wt F_R\times \wt F_R$ is a commutative division bisemiring.
\Ee
\vskip 11pt

\subsection{Bisemimodules}

\subsubsection{Definitions}

Let $R\Lr$ be a \lr semiring.  A {\bbf \lr $R\Lr$-semimodule $M\Lr$\/} is an additive abelian \lr semigroup (or monoid) $M\Lr$ together with a function $R_L\times M_L\to M_L$ \resp{$M_R\times  R\Rl\to M_R$} such that:
\Bean
\item $\begin{aligned}[t]
 r_L\ (g_{L_i}+g_{L_j})
&= r_L\ g_{L_i}+r_L\ g_{L_j}\;, \quad &&\forall\ r_L\in R_L\ , \; g_{L_i},g_{L_j}\in M_L\;, \\
\rresp{
(g_{R_i}+g_{R_j})\ r_R
&= g_{R_i}\ r_R+\ g_{R_j}\ r_R\;, \quad &&\forall\ r_R\in R_R\ , \; g_{R_i},g_{R_j}\in M_R}.\end{aligned}$

\item $\begin{aligned}[t]
 (r_L+s_L)\  g_{L_i}
 &=r_L\  g_{L_i}+s_L\  g_{L_i}\\
 \rresp{
  g_{R_i}\ (r_R+s_R)
 &=g_{R_i}\ r_R +  g_{R_i}\ s_R}\end{aligned}$
 
 \item $\begin{aligned}[t]
 r_L\ (s_L\ g_{L_i})
 &=(r_L\ s_L)\ g_{L_i}\\
 \rresp{
   ( g_{R_i}\ s_R)\ r_R
 &=g_{R_i}\ (s_R\ r_R)}\end{aligned}$
 \Ee
 
 {\bbf $M\Lr$ is a \lr unitary $R\Lr$-semimodule\/} if $R\Lr$ has an identity element $\ung\Lr$ such that
 \[ \ung_L\ g_{L_i}=g_{L_i}\;, \qquad
 \rresp{g_{R_i}\ung_R=g_{R_i}}.\]
 
 {\bbf $M\Lr$ is a \lr $R\Lr$-vector semispace\/} if:
 \Bi
 \item $R\Lr$ is a \lr division semiring.
 \item $M\Lr$ is a unitary \lr $R\Lr$-semimodule.
 \Ei
 \vskip 11pt
 
 \subsubsection{Definition}
 
 A {\bbf $R_R\times R_L$-bisemimodule $M\RL$\/} over a bisemiring $R\RL\equiv R_R\times R_L$ {\bbf is an additive abelian bisemigroup\/} (or bimonoid) $M\RL$ together with a (bi)function 
\begin{multline*}
(R\RL)\times (M\RL)\to M\RL\ \text{given by:}\\
  (r_R\times r_L)\times (g_{R_i}\times g_{L_i})
 \To (g_{R_i}\ r_R\times r_L\ g_{L_i})\;, \quad \forall\ (g_{R_i}\times g_{L_i})\in M\RL\;,\end{multline*}
 such that:
\Bean
\item $((g_{R_i}+g_{R_j})\ r_R\times 
r_L\ (g_{L_i}+g_{L_j}))
=(g_{R_i}\ r_R \times r_L\ g_{L_i})
\Times (g_{R_j}\ r_R \times r_L\ g_{L_j})$ where $\Times$ denotes the cross binary operation.

\item $(g_{R_i}\ r_R \times r_L\ g_{L_i})
\neq(g_{R_i}\ r_L \times r_R\ g_{L_i})$.
\Ee
\vskip 11pt

\subsubsection[Definition: tensor product of right and left semimodules over $R_R\times R_L$]{\bbf Definition: tensor product of right and left semimodules over $R_R\times R_L$}

The tensor product $M_R\otimes_{R_R\times R_L}M_L$ of the right $R_R$-semimodule $M_R$ by the left $R_L$-semimodule $M_L$ is the quotient $G\RL\big/ K\RL$ of the abelian bisemigroup $G\RL$ by the bisemisubgroup $K\RL$ generated by bielements of the form $(g_{R_i}+g_{R_j})\times (g_{L_i}+g_{L_j})$ submitted to the cross binary operation  ``$\Times$'' in such a way that:
\[ (g_{R_i}\times g_{L_i})
\Times (g_{R_j}\times g_{L_j})
\To
(g_{R_i}+ g_{R_j})
\times (g_{L_i}\times g_{L_j})\;.\]

\Bi
\item {\bf The coset\/} $(g_{R_i}\times g_{L_i})
+K\RL$ of the bielement
$(g_{R_i}\times g_{L_i})$ in $G\RL$ is noted {\bbf $(g_{R_i}\otimes g_{L_i})$} and {\bf is of the form\/}:
\[ \sum_\ell (n_\ell\ g_{R_{i_\ell}}
\otimes m_\ell\ g_{L_{i_\ell}})\;, \qquad
n_\ell,m_\ell\in \ZZ\;.\]

\item {\bbf The bigenerators\/} $(g_{R_i}\otimes g_{L_i})$ of $M_R\otimes_{R\RL}M_L$ satisfy the bilinearity condition:
\[
(g_{R_i} + g_{R_j})
\otimes (g_{L_i} + g_{L_j})
=(g_{R_i} \otimes g_{L_i})
+(g_{R_j} \otimes g_{L_j})
+(g_{R_i} \otimes g_{L_j})
+(g_{R_j} \otimes g_{L_i})\;.\]
\Ei
\vskip 11pt

\subsubsection{Proposition}

{\em $M_R\otimes_{R\RL}M_L$ is a $R_R\times R_L$-bisemimodule $M\RL$}
\vskip 11pt

\paragraph{Proof.} \quad Indeed, $M_R\otimes_{R\RL}M_L$  is a quotient bisemigroup $G\RL\big/ K\RL$ of which bielements are of the form $\sum_\ell (n_\ell\ g_{R_{i_\ell}}
\otimes m_\ell\ g_{L_{i_\ell}})$ according to definition 2.3.3.

On the other hand, a $R_R\times R_L$-bisemimodule $M\RL$ is an additive bisemigroup of which bielements are $(g_{R_i}\ r_R\times r_L\ g_{L_i})$ according to definition 2.3.2.

It then appears that $\sum_\ell (n_\ell\ g_{R_{i_\ell}}
\otimes m_\ell\ g_{L_{i_\ell}})$ is a general form for
$(g_{R_i}\ r_R\times r_L\ g_{L_i})$ if $r_R,r_L\in \ZZ$ and thus, we have the thesis.\epr
\vskip 11pt

\subsubsection{Proposition}

{\em Let $M_R$ \resp{$M_L$} be a \rl $R_R$-semimodule
\resp{$R_L$-semimodule}.

Let $A$ \resp{$B$} be a \rl $R$-semimodule over a ring $R$.

Then, $A\otimes_R B$ is a quotient semigroup $G/K$ while
$M_R\otimes_{R_R\times R_L}M_L$ is a quotient bisemigroup $G\RL\big/ K\RL$.}
\vskip 11pt

\paragraph{Proof.} \quad  Classically \cite{Hun}, \cite{Sch1},  \cite{Greu}, $A\otimes_R B$ is a quotient semigroup $G/K$ where $F$ is a free abelian semigroup on the set $A\times B$ and $K$ is a subsemigroup in such a way that the coset $(a,b)+K$ of the element $(a,b)$ in $F$ is denoted $a\otimes b$ and has the general form
\[ \sum_in_i\ (a_i\otimes b_i)\;, \qquad n_i\in\ZZ\ ,\; a_i\in A\ , \; b_i\in B\;.\]
So, $A\otimes_RB$ differs from $M_R\otimes_{R_R\times R_L}M_L$ in the sense that $A\otimes_R B$ has a linear character while $M_R\otimes_{R_R\times R_L}M_L$ behaves bilinearly.\epr
\vskip 11pt

\subsection{Inner product bisemispaces and bilinear Hilbert spaces}

\subsubsection{Definition: vector bisemispace}

\Bi
\item Recall that $M\Lr$ is a {\bbf \lr vector $R\Lr$-semispace\/} if:
\Bean
\item $R\Lr$ is a \lr division semiring.
\item $M\Lr$ is a unitary \lr $R\Lr$-semimodule.
\Ee

\item {\bbf A vector $R_R\times R_L$-bisemispace $M\RL$\/} is a unitary $R_R\times R_L$-bisemimodule over a division bisemiring $R_R\times R_L$.
\Ei
\vskip 11pt

\subsubsection{Proposition}

{\em The tensor product $M_R\otimes_{R_R\times R_L}M_L$ of a right vector $R_R$-semispace $M_R$ by a left vector $R_L$-semispace $M_L$ is a vector $R_R\times R_L$-bisemispace.}
\vskip 11pt

\paragraph{Proof.} \quad Indeed, $M_R\otimes_{R_R\times R_L}M_L$ is a unitary $ R_R\times R_L$-bisemimodule according to proposition 2.3.4.\epr
\vskip 11pt

\subsubsection{Proposition}

{\em {\bfseries The vector {\bbf $R_R\times R_L$\/}-bisemispace {\bbf $M_R\otimes_{R_R\times R_L}M_L$} of dimension {\bbf $n^2$} splits naturally into:
{\bbf\begin{align*}
M_R\otimes_{R_R\times R_L}M_L
&\equiv M_R\otimes M_L &&\text{(condensed form)}\\
&= (M_R\otimes_D M_L)
\oplus (M_R\otimes_{OD} M_L)\end{align*}}}
where:
\Bi
\item $(M_R\otimes_DM_L)$ is a {\bbf diagonal vector $R_R\times R_L$-bisemispace\/} of dimension $n$ characterized by a bilinear diagonal basis $\{e_\alpha \otimes f_\alpha \}_\alpha $, $\forall\ e_\alpha \in M_R$, $f_\alpha \in M_L$.

\item $(M_R\otimes_{OD}M_L)$ is an {\bbf off-diagonal vector $R_R\times R_L$-bisemispace\/} of dimension $(n^2-n)$ characterized by a bilinear off-diagonal basis $\{e_\alpha \otimes f_\beta  \}_{\alpha\neq\beta } $.
\Ei}
\vskip 11pt

\paragraph{Proof.} \quad The vector $R_R\times R_L$-bisemispace $M_R\otimes M_L$, being a unitary 
$R_R\times R_L$-bisemimodule, is an additive abelian bisemigroup.

So, its bielements $(x_{R_i}\times x_{L_i})$, $\forall\ x_{R_i}\in M_R$, $x_{L_i}\in M_L$, are submitted to the cross binary operation ``$\Times$'' according to definition 2.1.2.

Then, they split as follows:
\begin{multline*}
(x_{R_i}\times x_{L_i})
\Times (x_{R_j}\times x_{L_j})\\
\To
(x_{R_i}+ x_{R_j})\times (x_{L_i}+ x_{L_j})
=(x_{R_i}\times x_{L_i})+(x_{R_j}\times x_{L_j})
+(x_{R_i}\times x_{L_j})+(x_{R_j}\times x_{L_i})
\end{multline*}
in such a way that:
\Bi
\item $(x_{R_i}\times x_{L_i})$ and $(x_{R_j}\times x_{L_j})$ are diagonal bielements belonging to the diagonal vector $R_R\times R_L$-bisemispace $(M_R\otimes_DM_L)$ characterized by a bilinear diagonal basis $\{e_{\alpha_i} \otimes f_{\alpha_i} \}_{\alpha_i} $.

\item $(x_{R_i}\times x_{L_j})$ and $(x_{R_j}\times x_{L_i})$ are off-diagonal bielements belonging to the off-diagonal vector $R_R\times R_L$-bisemispace $(M_R\otimes_{OD}M_L)$ characterized by a bilinear off-diagonal basis $\{e_{\alpha_i} \otimes f_{\beta _i} \}_{\alpha_i\neq\beta _i} $.\epr
\Ei
\vskip 11pt

\subsubsection{Definitions: mixed product bisemispaces}

\Bean
\item A diagonal vector $R_R\times R_L$-bisemispace $(M_R\otimes_DM_L)$ endowed with a diagonal mixed product $[\centerdot,\centerdot]_D$ is called a {\bbf diagonal mixed product bisemispace\/} if the bilinear function
$[\centerdot,\centerdot]_D$ from $M_R\otimes_DM_L$ to
$R_R\times R_L$ satisfies the following conditions:
\Be
\item $[x_{R_i},x_{L_i}]_D\ge 0$
and $[x_{R_i},x_{L_i}]_D= 0$ if $ x_{R_i}=0$ or $x_{L_i}= 0$.

\item $[x_{R_i}+x_{R_j},x_{L_i}+x_{L_j}]_D
=[x_{R_i},x_{L_i}]_D +[x_{R_j},x_{R_j}]_D
+[x_{R_i},x_{L_j}]_D +[x_{R_j},x_{L_i}]_D$

\item $[x_{R_i},\alpha _{L_i}\ x_{L_i}]_D
=\alpha _{L_i}\ [x_{R_i},x_{L_i}]_D$, $\forall\ \alpha _{L_i}\in R_L$,

 $[\alpha _{R_i}\ x_{R_i},x_{L_i}]_D
=\alpha _{R_i}\ [x_{R_i},x_{L_i}]_D$, $\forall\ \alpha _{R_i}\in R_R$.
\Ee

The function $[\centerdot,\centerdot]_D$ will be called a diagonal mixed product.

\item Similarly, an  off-diagonal vector $R_R\times R_L$-bisemispace $(M_R\otimes_{OD}M_L)$ endowed with an  off-diagonal mixed product $[\centerdot,\centerdot]_{OD}$ will be called an {\bbf off-diagonal mixed product bisemispace\/} if the bilinear function
$[\centerdot,\centerdot]_{OD}$ from $M_R\otimes_{OD}M_L$ to
$R_R\times R_L$ satisfies conditions similar to those of a).
\Ee
\vskip 11pt

\subsubsection{Definition: extended mixed product bisemispace}

The vector $R_R\times R_L$-bisemispace $(M_R\otimes_{R_R\times R_L}M_L)$ (noted $M_R\otimes M_L$) endowed with an (extended) bilinear function $[\centerdot,\centerdot]$  will be called an extended mixed product bisemispace if:
\Bena
\item $M_R\otimes M_L$ is characterized by a bilinear basis $\{e^{\alpha} \otimes f_{\beta } \}_{\alpha,\beta} $, $\forall\ e^\alpha \in M_R$, $f_\beta \in M_L$.

\item the (extended) bilinear function $[\centerdot,\centerdot]$ from $M_R\times M_L$ to $R_R\times R_L$ satisfies conditions as those of definitions 2.4.4.
\Ee
\vskip 11pt

\subsubsection{Proposition (extended, diagonal and off-diagonal ``external'' \lr product bisemispaces)}

{\em Let $(M_R\otimes_{R_R\times R_L}M_L)$
(or $(M_R\otimes M_L)$), $(M_R\otimes_D M_L)$
and $(M_R\otimes_{OD} M_L)$ be respectively the extended, diagonal and off-diagonal mixed product bisemispaces endowed  with the corresponding bilinear functions 
$[\centerdot,\centerdot]$, 
$[\centerdot,\centerdot]_D$ and
$[\centerdot,\centerdot]_{OD}$.

Then, the {\bfseries linear projective morphisms\/}:
\begin{align*}
&\begin{array}{lllll}
p_L: & \qquad & M_R\otimes M_L  &\To & M_{R(P)/L}\equiv M_{R(P)}\otimes M_L\\
p_L: &\qquad & M_R\otimes_D M_L  &\To & M_{R(P)/_DL}\equiv M_{R(P)}\otimes_D M_L\\
p_L: & \qquad & M_R\otimes_{OD} M_L  &\To & M_{R(P)/_{OD}L}\equiv M_{R(P)}\otimes_{OD} M_L\end{array}\\
\noalign{\mbox{}\hfill where $p_L=\Hom_{R_R\times R_L}(M_R,M_L)$,}
\text{(resp. \quad}
&\begin{array}[t]{lllll}
p_R: & \qquad & M_R\otimes M_L  &\To & M_{L(P)/R}\equiv M_{L(P)}\otimes M_R\\
p_R: & \qquad & M_R\otimes_D M_L  &\To & M_{L(P)/_DR}\equiv M_{L(P)}\otimes_D M_R\\
p_R: & \qquad & M_R\otimes_{OD} M_L  &\To &  M_{L(P)/_{OD}R}\equiv M_{L(P)}\otimes_{OD} M_R\ )\end{array}\\
\noalign{\mbox{}\hfill where $p_R=\Hom_{R_R\times R_L}(M_L,M_R)$,}
\end{align*}
\vskip -15pt
{\bfseries transform\/} respectively {\bfseries the\/} extended, diagonal and off-diagonal {\bfseries mixed product bisemispaces into\/} extended, diagonal and off-diagonal {\bfseries external \lr product bisemispaces\/} $M_{R(P)/L}$, $M_{R(P)/_DL}$ and $M_{R(P)/{OD}L}$ \resp{$M_{L(P)/R}$, $M_{L(P)/_DR}$ and $M_{L(P)/{OD}R}$} by mapping the \rl vector $R_R$- \resp{$R_L$}-semispace $M_R$ \resp{$M_L$} onto the \lr vector $R_L$- \resp{$R_R$}-semispace $M_L$ \resp{$M_R$}, in such a way that $M_{R(P)}$ \resp{$M_{L(P)}$} becomes the dual semispace of $M_L$ \resp{$M_R}$.
}
\vskip 11pt

\paragraph{Proof.} \quad The projective linear morphisms 
\[
p_L=\Hom_{R_R\times R_L}(M_R,M_L)\qquad
\rresp{p_R =\Hom_{R_R\times R_L}(M_L,M_R)}\]
transform the bilinear functions $[\centerdot,\centerdot]$, $[\centerdot,\centerdot]_D$ and $[\centerdot,\centerdot]_{OD}$ respectively of $(M_R\otimes M_L)$, $(M_R\otimes_D M_L)$ and $(M_R\otimes_{OD} M_L)$ according to:

\[\begin{array}{llllll}
& p_L: & \qquad & [\centerdot,\centerdot]& \To & [\centerdot,\centerdot]^L\;, \\
& p_L: & \qquad & [\centerdot,\centerdot]_D& \To & [\centerdot,\centerdot]^L_D\;, \\
& p_L: & \qquad & [\centerdot,\centerdot]_{OD}& \To & [\centerdot,\centerdot]^L_{OD}\;, \\
\text{(resp.} \quad
& p_R: & \qquad & [\centerdot,\centerdot]& \To & [\centerdot,\centerdot]^R\;, \\
& p_R: & \qquad & [\centerdot,\centerdot]_D& \To & [\centerdot,\centerdot]^R_D\;, \\
& p_R: & \qquad & [\centerdot,\centerdot]_{OD}& \To & [\centerdot,\centerdot]^R_{OD}\;), \end{array}\]
where $[\centerdot,\centerdot]^L$, $[\centerdot,\centerdot]^L_D$ and $[\centerdot,\centerdot]^L_{OD}$ \resp{$[\centerdot,\centerdot]^R$, $[\centerdot,\centerdot]^R_D$ and $[\centerdot,\centerdot]^R_{OD}$} are respectively \lr extended, diagonal and off-diagonal external scalar products \cite{Sch2}:
\Bean
\item from $M_{R(P)}\times M_L$, $M_{R(P)}\times_D M_L$ and $M_{R(P)}\times_{OD} M_L$
\resp{$M_{R}\times M_{L(P)}$, $M_{R}\times_D M_{L(P)}$ and $M_{R}\times_{OD} M_{L(P)}$} respectively to $R_{R(P)}\times R_L$,  $R_{R(P)}\times_D R_L$ and $R_{R(P)}\times_{OD} R_L$
\resp{$R_{R}\times R_{L(P)}$, $R_{R}\times_D R_{L(P)}$ and $R_{R}\times_{OD} R_{L(P)}$} with:
\Bi
\item $R_{R(P)}$ \resp{$R_{L(P)}$} being the semiring
$R_{R}$ \resp{$R_{L}$} projected onto $R_{L}$ \resp{$R_{R}$}.
\item $M_{R(P)}$ \resp{$M_{L(P)}$} being $M_{R}$ \resp{$M_{L}$} projected onto $M_L$ \resp{$M_R$} in such a way that $M_{R(P)}$ \resp{$M_{L(P)}$} is the semispace of one-forms on 
$M_L$ \resp{$M_R$}.
\Ei

\item characterized by metrics of type $(1,1)$.
\Ee
$p_L$ \resp{$p_R$} is then a covariant functor.\epr
\vskip 11pt

\subsubsection{Proposition (extended, diagonal and off-diagonal ``internal'' \lr product bisemispaces)}

{\em Let $M_{R(P)/L}$
$M_{R(P)/_DL}$
and $M_{R(P)/_{OD}L}$ 
\resp{$M_{L(P)/R}$
$M_{L(P)/_DR}$
and $M_{L(P)/_{OD}R}$}
be respectively the extended, diagonal and off-diagonal external \lr product bisemispaces.

Then, {\bfseries the bijective linear isometric homomorphisms\/}:
\[ \mbox{\bbf$B_L$}: \quad M_{R(P)}\To M_{L_R}
\rresp{ \mbox{\bbf$B_R$}: \quad M_{L(P)}\To M_{R_L}},\]
{\bfseries mapping each covariant (resp. contravariant) element of 
$M_{R(P)}$ \resp{$M_{L(P)}$} onto the corresponding contravariant \resp{covariant} element\/} of $M_L$ \resp{$M_R$}, transform these external \lr product bisemispaces into the corresponding internal \lr product bisemispaces according to:
\begin{align*}
&\begin{array}{lllll}
B_L: & \qquad & M_{R(P)/L}  &\To & M_{L_R/L}\equiv M_{L_R}\otimes M_L\\
B_L: & \qquad & M_{R(P)/_DL}  &\To & M_{L_R/_DL}\equiv M_{L_R}\otimes_D M_L\\
B_L: & \qquad & M_{R(P)/_{OD}L}  &\To & M_{L_R/_{OD}L}\equiv M_{L_R}\otimes_{OD} M_L\\\end{array}\\
\text{(resp. \quad}
&\begin{array}[t]{lllll}
B_R: & \qquad & M_{L(P)/R}  &\To & M_{R_L/R}\equiv M_{R_L}\otimes M_R\\
B_R: & \qquad & M_{L(P)/_DR}  &\To & M_{R_L/_DR}\equiv M_{R_L}\otimes_D M_R\\
B_R: & \qquad & M_{L(P)/_{OD}R}  &\To & M_{R_L/_{OD}R}\equiv M_{R_L}\otimes_{OD} M_R
\ ).\end{array}
\end{align*}
}
\vskip 11pt

\paragraph{Proof.} \quad Indeed, the bijective linear isometric homomorphisms $B_L$ \resp{$B_R$} transform the external scalar products of the external \lr product bisemispaces 
$M_{R(P)/L}$, $M_{R(P)/_DL}$ and $M_{R(P)/_{OD}L}$
\resp{$M_{L(P)/R}$, $M_{L(P)/_DR}$ and $M_{L(P)/_{OD}R}$} according to:
\[\begin{array}{llllll}
& B_L: & \qquad & [\centerdot,\centerdot]^L& \To &  \langle\centerdot,\centerdot\rangle^L\;, \\
& B_L: & \qquad & [\centerdot,\centerdot]^L_D& \To &  \langle\centerdot,\centerdot\rangle^L_D\;, \\
& B_L: & \qquad & [\centerdot,\centerdot]^L_{OD}& \To &  \langle\centerdot,\centerdot\rangle^L_{OD}\;, \\
\text{(resp. \quad}
& B_R: & \qquad & [\centerdot,\centerdot]^R& \To &  \langle\centerdot,\centerdot\rangle^R\;, \\
& B_R: & \qquad & [\centerdot,\centerdot]^R_{D}& \To &  \langle\centerdot,\centerdot\rangle^R_{D}\;, \\
& B_R: & \qquad & [\centerdot,\centerdot]^R_{OD}& \To &  \langle\centerdot,\centerdot\rangle^R_{OD}\;), \end{array}\]
where $\langle\centerdot,\centerdot\rangle^L$,
  $\langle\centerdot,\centerdot\rangle^L_D$
and  $\langle\centerdot,\centerdot\rangle^L_{OD}$
\resp{$\langle\centerdot,\centerdot\rangle^R$,
 $\langle\centerdot,\centerdot\rangle^R_{D}$ and
 $\langle\centerdot,\centerdot\rangle^R_{OD}$} are respectively \lr extended, diagonal and off-diagonal inner (or internal) scalar products
from $M_{L_R}\times M_L$,
$M_{L_R}\times_D M_L$ and
$M_{L_R}\times_{OD} M_L$
\resp{$M_{R_L}\times M_R$,
$M_{R_L}\times_D M_R$ and
$M_{R_L}\times_{OD} M_R$}
respectively to $R_{L_R}\times R_L$, $R_{L_R}\times_D R_L$,
and $R_{L_R}\times_{OD} R_L$
\resp{$R_{R_L}\times R_R$, $R_{R_L}\times_D R_R$
and $R_{R_L}\times_{OD} R_R$} with
$R_{L_R}$ \resp{$R_{R_L}$} being the semiring $R_R$ \resp{$R_L$} projected onto $R_L$ \resp{$R_R$} under the condition that the metrics of these inner scalar products \cite{Wei} be of type $(0,2)$ \resp{$(2,0)$}.
\vskip 11pt

So, $B_L$ \resp{$B_R$}, corresponding to the Riesz lemma mapping \cite{G-S}, \cite{R-N}, transforms the metrics of type $(1,1)$ of the external scalar products into metrics of type $(0,2)$ \resp{$2,0)$} of the respective inner scalar products \cite{B-R}, \cite{C-H}, \cite{God}, \cite{J-N}.\epr
\vskip 11pt

\subsubsection{Corollary}

{\em Every  \lr diagonal internal product bisemispace $(M_{L_R}\otimes_D M_L)$ \resp{$(M_{R_L}\otimes_D M_R)$}, endowed with the diagonal inner scalar product $\langle \centerdot,\centerdot\rangle^L_D$ \resp{$\langle \centerdot,\centerdot\rangle^R_D$}, is a normed bilinear semispace with ``bilinear norm''
$\langle x_{L_{R_i}},x_{L_i}\rangle^L_D$
\resp{$\langle x_{R_{L_i}},x_{R_i}\rangle^R_D$}, $x_{L_{R_i}}\in M_{L_R}$,
$x_{L_{i}}\in M_{L}$,
$x_{R_{L_i}}\in M_{R_L}$,
$x_{R_{i}}\in M_{R}$ in such a way that:
\Bena
\item the semispaces $M_L$, $M_{L_R}$, $M_R$ and $M_{R_L}$ are normed linear semispaces with norm 
\[\|x_{L_i}\| = (\langle x_{L_{R_i}},x_{L_i}\rangle^L_D)^\half\;.\]
\item $\|x_{L_i}\| = \|x_{L_{R_i}}\|=\|x_{R_{L_i}}\|=\|x_{R_i}\|$.
\Ee}
\vskip 11pt

\paragraph{Proof.} \quad
\Bena
\item The requirement that the internal product bisemispaces
$(M_{L_R}\otimes_D M_L)$ and $(M_{R_L}\otimes_D M_R)$ be normed bilinear semispaces is similar to the condition that every inner product space is a normed linear space with norm $\|x\|=\langle x,x\rangle$ where $\langle\centerdot,\centerdot\rangle$ denotes a classical inner scalar product.

\item The equalities $\|x_{L_i}\| = \|x_{L_{R_i}}\|=\|x_{R_{L_i}}\|=\|x_{R_i}\|$ result  from the fact that these elements, being symmetric by construction, have a same length.

\item As $(M_{L_R}\otimes_D M_L)$ and $(M_{R_L}\otimes_D M_R)$ are normed bilinear semispaces, we have that:
\Be
\item if their elements $x_{L_{R_i}}$, $x_{L_{i}}$, $x_{R_{L_i}}$ and $x_{R_{i}}$ are real-valued, then the semirings $R_L$ and $R_{R_L}$ are semirings of positive real numbers while the semirings $R_R$ and $R_{R_L}$ are semirings of negative real numbers (see proof of proposition 2.4.7 for the definition of these semirings in the present context).

\item if their elements $x_{L_{R_i}}$, $x_{L_{i}}$, $x_{R_{L_i}}$ and $x_{R_{i}}$ are complex-valued, then:
\Be
\item[$\alpha $)] the semirings $R_L$ and $R_{R_L}$ are semirings of complex numbers $z=a+ib$;
\item[$\beta $)] the semirings $R_{L_R}$ and $R_{R}$ are semirings of conjugate complex numbers\linebreak $z^*=a-ib$;
\item[$\gamma $)] the \lr diagonal inner scalar product
$\langle x_{L_{R_i}},x_{L_i}\rangle^L_D$
\resp{$\langle x_{R_{L_i}},x_{R_i}\rangle^R_D$} is such that
\Bi
\item $ x_{L_{R_i}}=\o{x_{L_i}}$
\resp{$ x_{R_{L_i}}=\o{x_{R_i}}$}

\item $ \langle y_{L_{R_i}},{x_{L_i}}\rangle^L_D
= \langle {x_{R_{L_i}},{y_{R_i}}}\rangle^R_D$, $\forall\ 
y_{L_{R_i}}\in M_{L_R}$, $y_{R_i}\in M_R$.
\Ei
\Ee 
Thus, $\langle x_{L_{R_i}},x_{L_i}\rangle^L_D$ is a sesquilinear form, linear on the right and antilinear on the left.\epr
\Ee
\Ee
\vskip 11pt

\subsubsection{Proposition}

{\em The diagonal and extended internal \lr product bisemispaces $M_{L_R/_DL}$ (resp. \linebreak {$M_{R_L/_DR}$)}
and $M_{L_R/L}$ \resp{$M_{R_L/R}$} are respectively {\bfseries an orthogonal and an extended \lr separable bilinear Hilbert semispace\/} while the off-diagonal internal \lr product bisemispace $M_{L_R/_{OD}L}$ \resp{$M_{R_L/_{OD}R}$} has been called an (internal) \lr bilinear magnetic space.
}
\vskip 11pt

\paragraph{Proof.}\quad
\Bean  
\item 
\Bi 
\item The diagonal internal \lr product bisemispace
$M_{L_R/_DL}\equiv M_{L_R}\otimes_D M_L$ \resp{$M_{R_L/_DR}\equiv M_{R_L}\otimes_D M_R$}, endowed with a diagonal \lr inner scalar product $\langle\centerdot,\centerdot\rangle^L_D$
\resp{$\langle\centerdot,\centerdot\rangle^R_D$} verifying:
\[ \langle\centerdot,\centerdot\rangle^R_D
= \langle{\o{\centerdot,\centerdot}}\rangle^L_D\]
under a (bi)involution, is characterized by:
\Be
\item bielements 
\begin{align*}
(x_{L_{R_i}}\otimes_D x_{L_i})
&=\sum\limits^\ell_{\alpha =1} (x_{L_{R_{i_\alpha }}}\ e_\alpha \otimes_D x_{L_{i_\alpha }}\ e_\alpha) \\
&=\sum\limits^\ell_{\alpha =1} x_{L_{R_{i_\alpha }}}\centerdot x_{L_{i_\alpha }}\ (e_\alpha \otimes e_\alpha )\;,
&& 1\le\alpha \le\ell\le\infty \;, \end{align*}
developed in the bilinear orthogonal basis $\{e_\alpha \otimes e_\alpha\}_{\alpha =1}^\ell$ of dimension $\ell$.

The bielements $(x_{L_{R_i}}\otimes_D x_{L_i})$
are two confounded elements.

\item left inner scalar products $\langle x_{L_{R_i}}, x_{L_i}\rangle^L_D$.

If
\[ \o x_{L_{R_i}}\centerdot x_{L_i}
= \{\o x_{L_{R_{i_1}}}\centerdot x_{L_{i_1}},\dots,
\o x_{L_{R_{i_\alpha }}}\centerdot x_{L_{i_\alpha }},\dots,
\o x_{L_{R_{i_\ell }}}\centerdot x_{L_{i_\ell }}\}\]
denotes the set of $\ell$-bituple of complex numbers in $\CC^\ell\times_D \CC^\ell$, then $\langle x_{L_{R_i}}, x_{L_i}\rangle^L_D$ can be developed according to:
\[ \langle x_{L_{R_i}}, x_{L_i}\rangle^L_D=
\sum^\ell_{\alpha =1} \o x_{L_{R_{i_\alpha }}}\centerdot x_{L_{i_\alpha }}\]
(see corollary 2.4.8).
\Ee

\item {\bbf The diagonal internal left product bisemispace 
$M_{L_R/_DL}\equiv M_{L_R}\otimes_D M_L$ is in one-to-one correspondence with the linear inner product semispace $\Hs$ which is a pre-Hilbert semispace\/} because:
\Be
\item $\Hs$ is also characterized by inner scalar products, noted generally $\langle\centerdot,\centerdot\rangle$, from $\Hs\times \Hs$ to $\CC$.  These inner scalar products of  $\Hs$ are the left inner scalar products
$\langle\centerdot,\centerdot\rangle^L_D$ of $M_{L_R/_DL}$.

\item the bielements $(x_{L_{R_i}}\otimes_D x_{L_i})$ of $M_{L_R/_DL}$ are in one-to-one correspondence with the elements $x_{L_i}$ of $\Hs$.
\Ee

So, it can be understood why the diagonal internal left product bisemispace $M_{L_R/_DL}$ is called (and is) an orthogonal  left bilinear pre-Hilbert semispace (it is an orthogonal left bilinear Hilbert semispace if it is complete with respect to $\langle\centerdot,\centerdot\rangle^L_D$).

\item Similarly, the diagonal internal right product bisemispace 
$M_{R_L/_DR}\equiv M_{R_L}\otimes_D M_R$ is an orthogonal right bilinear Hilbert semispace because it would have to be in one-to-one correspondence with a linear inner product semispace
$\Hs_R$ (which would be a right Hilbert semispace) characterized by:
\Be
\item elements $x_{{R_i}}$  in one-to-one correspondence with the bielements $(x_{R_{L_i}}\otimes_Dx_{R_i})$ of $M_{R_L/_DR}$ (in fact, the elements $x_{L_{R_i}}$ can be identified with the elements $x_{R_{i}}$).

\item inner scalar products $\langle\o{\centerdot,\centerdot}\rangle$ from $\Hs_R\times \Hs_R$ to $\CC$.
\Ee

Then, it becomes clear that $\Hs_R\equiv \Hs$ if the elements $x_{{R_i}}$ of $\Hs_R$ are identified with the elements $x_{{L_i}}$ of $\Hs$.
\Ei  

\item The extended internal \lr product bisemispace
$M_{L_R/L}\equiv M_{L_R}\otimes M_L$
\resp{$M_{R_L/R}\equiv M_{R_L}\otimes M_R$}, endowed with an extended \lr inner product 
$\langle{\centerdot,\centerdot}\rangle^L$
\resp{$\langle{\centerdot,\centerdot}\rangle^R$} verifying
$\langle {\centerdot,\centerdot}\rangle^R=
\langle {\centerdot,\centerdot}\rangle^L$, is characterized by:
\Be
\item bielements 
\begin{align*}
(x_{L_{R_i}}\otimes x_{L_i})
&=\sum^\ell_{\alpha =1} \sum^\ell_{\beta  =1}
(x_{L_{R_{i_\alpha }}}\ e_\alpha \otimes x_{L_{i_\beta  }}\ e_\beta ) \\
&=\sum^\ell_{\alpha =1} \sum^\ell_{\beta  =1}
x_{L_{R_{i_\alpha }}}\centerdot x_{L_{i_\beta  }}\ (e_\alpha \otimes e_\beta  )\;,
&& 1\le\alpha,\beta  \le\ell\le\infty \;, \end{align*}
developed in the bilinear basis 
$\{e_\alpha \otimes e_\beta \}^\ell_{\alpha ,\beta =1}$.

\item left extended inner scalar products
$\langle x_{L_{R_i}}, x_{L_i}\rangle^L$ which can be developed according to:
\[ \langle x_{L_{R_i}}, x_{L_i}\rangle^L
= \sum^\ell_{\alpha =1} \sum^\ell_{\beta  =1}
\o x_{L_{R_{i_\alpha }}}\centerdot x_{L_{i_\beta }}\]
if $\o x_{L_{R_{i }}}\centerdot x_{L_{i }}$ denotes the set of $\ell^2$-bituples of $\CC^\ell\times\CC^\ell$.

With reference to a), it is then clear that $M_{L_R/L}$
\resp{$M_{R_L/R}$} is called an extended \lr bilinear Hilbert semispace since it is endowed with an extended \lr inner scalar product.
\Ee

\item Similarly, the off-diagonal internal \lr product bisemispace 
$M_{L_R/_{OD}L}$ \resp{$M_{R_L/_{OD}R}$}
 is called an internal \lr bilinear magnetic semispace because it is endowed with:
 \Be
 \item a \lr off-diagonal inner scalar product
$ \langle{\centerdot,\centerdot}\rangle^L_{OD}$
\resp{$ \langle{\centerdot,\centerdot}\rangle^R_{OD}$}
developed according to:
\[ \langle x_{L_{R_i}}, x_{L_i}\rangle^L_{OD}
= \underset{\alpha \neq \beta }{\sum^\ell_{\alpha =1} \sum^\ell_{\beta  =1}}
\o x_{L_{R_{i_\alpha }}}\centerdot x_{L_{i_\beta }}\]
if $\o x_{L_{R_{i }}}\centerdot x_{L_{i }}$ denotes the set of bituples of $\CC^\ell\times_{OD}\CC^\ell$.

\item off-diagonal bielements
\begin{align*}
(x_{L_{R_i}}\otimes_{OD} x_{L_i})
&=\underset{\alpha \neq\beta }{\sum^\ell_{\alpha =1} \sum^\ell_{\beta  =1}}
x_{L_{R_{i_\alpha }}}\centerdot x_{L_{i_\beta  }}\ (e_\alpha \otimes_{OD} e_\beta ) \;,
 \end{align*}
 in the bilinear non-orthogonal basis 
$\{e_\alpha \otimes_{OD} e_\beta \}_{\alpha \neq\beta }$.
\epr
\Ee
\Ee
\vskip 11pt

\subsection{Semialgebras and bisemialgebras}

\subsubsection{Definitions}

\Be
\item {\bbf Left and right semialgebras}:

Let $R_L$ \resp{$R_R$} be a commutative \lr semiring with identity.  A \lr $R_L$-semialgebra $A_L$
\resp{$R_R$-semialgebra $A_R$} is a semiring $A_L$ \resp{$A_R$} such that:
\Be
\item $(A_L,+)$ \resp{$(A_R,+)$} is a unitary \lr $R_L$-semimodule \resp{$R_R$-semimodule}.

\item $\mu _L:A_L\otimes A_L \to A_L$
\resp{$\mu _R:A_R\otimes A_R \to A_R$} is a linear homomorphism verifying:
\begin{align*}
r_L\ (a_L\ b_L)=(r_L\ a_L)\ b_L=a_L\ (r_L\ b_L)\;, &&\forall\ r_L\in R_L\ , \; a_L,b_L\in A_L\;,\\
\rresp{
 (a_R\ b_R)\ r_R=a_R\ (b_R\ r_R)=b_R\ (a_R\ r_R)\;, &&\forall\ r_R\in R_R\ , \; a_R,b_R\in A_R}\;.\end{align*}

\item $\eta _L:R_L\to A_L$
\resp{$\eta _R:R_R\to A_R$} is an injective homomorphism.
\Ee

$A_R$ is thus the opposite semialgebra of $A_L$.

A \lr $R_L$-semialgebra $A_L$ \resp{$R_R$-semialgebra $A_R$} which, as a  \lr semiring, is a \lr division semiring, is called a \lr division semialgebra.

\item Let $R\RL=R_R\times R_L$ be a bisemiring as introduced in definition 2.2.1.

Then, a {\bbf $R\RL$-bisemialgebra $A_R\otimes A_L$\/} is a bisemiring $A\RL$ such that:
\Be
\item $(A\RL,+)$ is a unitary $R\RL$-bisemimodule;
\item $\mu \RL:A\RL\Times A\RL\to A\RL$ is a bilinear homomorphism such that $\Times$ is the cross binary operation acting on the bielements $(a_R\times a_L)\in A\RL$ as follows:
\[
[(a_R\times a_L)\Times (b_R\times b_L)]
=[(a_R+b_R)\times  (a_L+ b_L)]\;.\]

\item $\eta \RL:R\RL\to A\RL$ is an injective homomorphism.
\Ee
\Ee
\vskip 11pt

\subsubsection{Definition: left and right cosemialgebras}

Let $A_L$ \resp{$A_R$} be a \lr $R_L$-semialgebra
\resp{$R_R$-semialgebra}.  It is also a unitary \lr 
$R_L$-semimodule
\resp{$R_R$-semimodule}.

Let $p_L$ \resp{$p_R$} be the linear projective morphism:
\[ p_L: \quad A_R\To A_{R(P)}\;, 
\rresp{p_R: \quad A_L\To A_{L(P)}},\]
mapping $A_R$ \resp{$A_L$} into
$A_{R(P)}$ \resp{$A_{L(P)}$} onto
$A_L$ \resp{$A_R$} as introduced in proposition 2.4.6.

By this way, the $R(P)$- \resp{$L(P)$}-semialgebra 
$A_{R(P)}$ \resp{$A_{L(P)}$} becomes the dual semialgebra or {\bbf cosemialgebra of $A_L$ \resp{$A_R$}\/} such that:
\Bean
\item $
 \Delta _L: \quad A_{R(P)}\to A_{R(P)}\times A_L$
\quad  \resp{$\Delta _R: \quad A_{L(P)}\to A_{L(P)}\times A_R$}

is a linear homomorphism called comultiplication.

\item $\varepsilon _L: \quad A_{R(P)}\to R(P)$
\quad  \resp{ $\varepsilon _R: \quad A_{L(P)}\to L(P)$} 

is a linear form.
\Ee
\vskip 11pt

\subsubsection[Definition: bisemialgebras $(A_{R(P)}\otimes A_L)$ and $(A_{L(P)}\otimes A_R)$]{\bbf Definition: bisemialgebras $(A_{R(P)}\otimes A_L)$ and $(A_{L(P)}\otimes A_R)$}

Let  $R_{R(P)\times L}$ \resp{$R_{L(P)\times R}$} be the bisemiring obtained by projecting the semiring $R_R$ \resp{$R_L$} onto $R_L$ \resp{$R_R$}.

Then, the $R_{R(P)\times L}$-bisemialgebra $A_{R(P)}\otimes A_L$
\resp{$R_{L(P)\times R}$-bisemialgebra $A_{L(P)}\otimes A_R$} is a bisemiring
$A_{R(P)\times L}$
\resp{$A_{L(P)\times R}$} in such a way that:
\Bean
\item $(A_{R(P)\times L},+)$
\resp{$(A_{L(P)\times R},+)$} is a unitary
$R_{R(P)\times L}$- \resp{$R_{L(P)\times R}$}-bisemimodule.

\item $\begin{aligned}[t]
&\mu _{R(P)\times L}: \quad A_{R(P)\times L}\Times
A_{R(P)\times L}\to A_{R(P)\times L}\\
\text{(resp.} \quad
&\mu _{L(P)\times R}: \quad A_{L(P)\times R}\Times
A_{L(P)\times R}\to A_{L(P)\times R}\ )\end{aligned}$

is a bilinear homomorphism where $\Times$ is the cross binary operation.

\item $\eta _{R(P)\times L}:R_{R(P)\times L}\to A_{R(P)\times L}$
\quad  \resp{$\eta _{L(P)\times R}:R_{L(P)\times R}\to A_{L(P)\times R}$} 

is an injective homomorphism.
\Ee
\vskip 11pt

\subsubsection{Proposition}

{\em {\bfseries A bisemialgebra of Hopf is given by the triple
$((A_{R(P)}\otimes A_L),S_b\ (A_{L(P)}\otimes A_R))$\/} where:
\Bi
\item $(A_{R(P)}\otimes A_L)$ and $(A_{L(P)}\otimes A_R)$ are respectively the left and right bisemialgebras introduced in definition 2.5.3.

\item $S_b$ \resp{$S_b^{-1}$} is a bilinear antipode
\[ S_b: \quad A_{R(P)}\otimes A_L\to A_{L(P)}\otimes A_R
\rresp{S^{-1}_b: \quad A_{L(P)}\otimes A_R\to A_{R(P)}\otimes A_L}\]
mapping bijectively the \lr bisemialgebra into its symmetric \rl equivalent.
\Ei
}
\vskip 11pt

\paragraph{Proof.} \quad
\Be
\item Recall that {\bbf a Hopf algebra\/} over a field $k$ is a bialgebra $H(\mu ,\eta ,\Delta ,\varepsilon )$ in such a way that:
\Be
\item $H(\mu ,\eta )$ is a $k$-algebra where $\mu =H\otimes H \to H$ is linear and $\eta :k\to H$ is an injective homomorphism.

\item $H(\Delta ,\varepsilon )$ is a coalgebra in such a way that $\Delta :H\to H\otimes H$ is a linear morphism called comultiplication and $\varepsilon =H\to k$ is a linear form.

\item There exists a linear map $S:H\to H$ called antipode based on the convolution product $*$ on $\Hom_k(H,H)$ defined by \cite{Cae}, \cite{M-M}, \cite{Car}:
\[\forall h_1,h_2: \quad H\to H \qquad
h_1*h_2 = \mu \circ (h_1\otimes h_2)\circ\Delta \;.\]
\Ee

\item The corresponding bilinear case is naturally a bisemialgebra of Hopf where $A_{R(P)}$ \resp{$A_{L(P)}$} is the \lr cosemialgebra over $A_L$ \resp{$A_R$}.

Consequently, the bilinear antipode $S_b$ \resp{$S_b^{-1}$}, being an antihomomorphism of bisemialgebra of Hopf, maps naturally the \lr bisemialgebra of Hopf 
$A_{R(P)}\otimes A_L$
\resp{$A_{L(P)}\otimes A_R$} into the symmetric \rl bisemialgebra of Hopf
$A_{L(P)}\otimes A_R$
\resp{$A_{R(P)}\otimes A_L$}.\epr
\Ee
\vskip 11pt

\subsubsection{Definition: $*$-bisemialgebra}

Let $R_L$ \resp{$R_R$} denote the field of complex \resp{conjugate complex} numbers.  Then, a $R_R\times R_L$-bisemialgebra $A_R\otimes A_L$ is a $*$-bisemialgebra if there exist:
\Bean
\item the right involution:
\[ I_{L\to R}: \quad A_L\to A_R\;, \quad a_L\to a^*_R\;, \quad \forall\ a_L\in A_L\ , \; a^*_L\in A_R\;.\]

\item the left involution:
\[ I_{R\to L}: \quad A_R\to A_L\;, \quad a^*_R\to a_L\;, \]
\Ee
verifying $I_{R\to L}=I^{-1}_{L\to R}$ such that $a^*_L=(a^*_R)^*$.

So, the involution maps \cite{Dix} are clearly defined in $*$-bisemialgebras.

\section{Bilinear algebraic semigroups and bisemischemes}

The fundamental \SMST and \BSMST, having been introduced in chapter 2, it is now the aim of chapter 3 to take up some concrete cases of \BSMST, under the circumstances especially bilinear algebraic semigroups (as announced in chapter 1) and bisemischemes.
\vskip 11pt

\subsection{Bilinear algebraic semigroups}

\subsubsection[The Gauss decomposition of the $R-R$-bimodule $\GL_n(R)$]{\bbf The Gauss decomposition of the $R-R$-bimodule $\GL_n(R)$}

Let $R$ be a ring.  The set $M_n(R)$ of all $(n\times n)$ square matrices over $R$ is a left $R$-module under addition of matrices.   But, $M_n(R)$ is also a $R-R$-bimodule under addition of matrices because every $(n\times n)$ square matrix of $M_n(R)$ is the product of a $(n\times 1)$ column matrix by a $(1,n)$ row matrix: so, $M_n(R)$ is the product of the left $R$-module of column matrices by the right $R$-module of row matrices.

Similarly, the group $\GL_n(R)$ of all $(n\times n)$ invertible matrices $\gl_n(R)$ is also a $R-R$-bimodule.

Furthermore, the group $\GL_n(R)$ has the Gauss linear  decomposition for every regular matrix $\gl_n(R)$:
\[ \gl_n(R)=\delta (R)\ \xi _R(R)\ \xi _L(R)\]
where:
\Bi
\item $\delta (R)$ is a diagonal matrix of dimension $n$.
\item $\xi _R (R)$ is a lower unitriangular matrix.
\item $\xi _L (R)$ is an upper unitriangular matrix.
\Ei
\vskip 11pt

\subsubsection[The algebraic bilinear semigroup of matrices $\GL_n(\wt F_R\times \wt F_L)$]{\bbf The algebraic bilinear semigroup of matrices $\GL_n(\wt F_R\times \wt F_L)$}

Taking into account the general existence of semiobjects and bisemiobjects, we have to consider:
\Bean
\item instead of a ring $R$, {\bbf a triple $(R_L,R_R,R\RL)$} where:
\Bi
\item $R_L$ and $R_R$ are respectively left and right semirings.
\item $R\RL =R_R\times R_L$ is the associated bisemiring.
\Ei

\item that $\GL_n(R)$, being a $R-R$-bimodule, becomes over $R_R\times R_L$ a $R\RL$-bisemimodule $\GL_n(R_L\times R_L)$, introduced in definition 2.3.2.
\Ee

More concretely let $\wt F_R$ and $\wt F_L$ be the right and left algebraic symmetric extensions of a global number field $k$ of characteristic $0$, as introduced in section 1.2.  Then, the algebraic group with entries in
$(\wt F_R\times \wt F_L)$ is a bilinear algebraic semigroup written according to $\GL_n(\wt F_R\times \wt F_L)$. Its elements are $(n\times n)$ square regular matrices forming a $\wt F_R\times \wt F_L$-bisemimodule under addition.
\vskip 11pt

As in section 3.1.1, the bilinear algebraic  semigroup 
$\GL_n(\wt F_R\times \wt F_L)$ can be decomposed into the product of the right semigroup $T^t_n(\wt F_R)$ of lower triangular matrices with entries in $\wt F_R$ by the left semigroup 
$T_n(\wt F_L)$ of upper triangular matrices with entries in $\wt F_L$ according to
\[ \GL_n(\wt F_R\times \wt F_L)
= T^t_n(\wt F_R)\times T_n(\wt F_L)\;.\]
$T_n(\wt F_L)$ \resp{$T^t_n(\wt F_R)$} can be viewed as an operator
\[ T_n(\wt F_L): \quad\wt F_L\To T^{(n)}(\wt F_L)\equiv V_L\qquad
\rresp{T^t_n(\wt F_R): \quad\wt F_R\To T^{(n)}(\wt F_R)\equiv V_R}\]
sending the \lr semifield $\wt F_L$ \resp{$\wt F_R$} into the \lr 
$T_n(\wt F_L)$-semimodule $T^{(n)}(\wt F_L)$
\resp{$T^t_n(\wt F_R)$-semimodule $T^{(n)}(\wt F_R)$}
 which is a \lr affine semispace $V_L$ \resp{$V_R$} of dimension $n$ restricted to the upper \resp{lower} half space.
 \vskip 11pt
 
 \subsubsection{Proposition}
 
 {\em The algebraic bilinear semigroup of matrices 
 $\GL_n(\wt F_R\times \wt F_L)$ can be viewed as an operator:
 \[\GL_n(\wt F_R\times \wt F_L): \qquad
\wt F_R\times \wt F_L\To
G^{(n)}(\wt F_R\times \wt F_L)\equiv V_R\otimes_{\wt F_R\times \wt F_L}V_L\]
sending the bisemifield $\wt F_R\times \wt F_L$ into the
$\GL_n(\wt F_R\times \wt F_L)$-bisemimodule
$G^{(n)}(\wt F_R\times \wt F_L)$ which is:
\Bean
\item the representation space of the algebraic bilinear semigroup of matrices $\GL_n(\wt F_R\times \wt F_L)$;
\item a $n^2$-dimensional affine bisemispace $(V_R\otimes_{\wt F_R\times \wt F_L}V_L)$.
\Ee}
\vskip 11pt

\paragraph{Proof.} \quad
\Be
\item The $n^2$-dimensional affine bisemispace $G^{(n)}(\wt F_R\times \wt F_L)$ is naturally the tensor product over $(\wt F_R\times \wt F_L)$ of the right $n$-dimensional affine semispace 
$T^{(n)}(\wt F_R)$ by its left symmetric equivalent $T^{(n)}(\wt F_L)$ in such a way that its bielements are general algebraic bipoints $(-a_1\times_D a_1,\dots,
-a_n\times_D a_n)$ which are $n^2$-bituples as developed in chapter 1.

\item In fact, $\GL_n(\wt F_R\times \wt F_L)$ acts into the
$\GL_n(\wt F_R\times \wt F_L)$-bisemimodule
$G^{(n)}(\wt F_R\times \wt F_L)$ on an irreducible $n^2$-dimensional affine bisemispace 
$P^{(n)}(\wt F_R\times \wt F_L)$ considered as the (irreducible) unitary representation space of
$\GL_n(\wt F_R\times \wt F_L)$ \cite{Pie1} and as its (bi)center.\epr
\Ee
\vskip 11pt

\subsubsection{Proposition}

{\em {\bfseries The bilinear algebraic semigroup {\bbf $\GL_n(\wt F_R\times \wt F_L)$} has the following Gauss bilinear decomposition\/}:
{\bbf \[
 \GL_n(\wt F_R\times \wt F_L)
 =[D_n(\wt F_R)\times D_n(\wt F_L)]
 \times [UT^t_n(\wt F_R)\times UT_n(\wt F_L)]\]}
 where:\Bi
 \item $D_n(\centerdot)$ is the subgroup of diagonal matrices;
 \item $UT_n(\centerdot)$ \resp{$UT^t_n(\centerdot)$} is the subgroup of upper \resp{lower} unitriangular matrices.
 \Ei
 }
 \vskip 11pt
 
 \paragraph{Proof.} \quad This follows directly from the decomposition
 \[ \GL_n(\wt F_R\times \wt F_L)
 =T^t_n(\wt F_R)\times T_n(\wt F_L)\]
 of $\GL_n(\wt F_R\times \wt F_L)$ into the product of the semigroup $T^t_n(\wt F_R)$ of lower triangular matrices by the semigroup $T_n(\wt F_L)$ of upper triangular matrices.\epr
 \vskip 11pt

 \subsubsection{Proposition}
 
 {\em 
 Let $F=F_R\cup F_L$ denote the set of completions associated with the finite algebraic (closed) extension 
 $\wt F=\wt F_R\cup \wt F_L$ of the number field $k$.
 
 Let $\gl_n(F)=[u^t_n(F)\times u_n(F)]\times d_n(F)$ be the Gauss (linear) decomposition of the matrix 
 $\gl_n(F)$ of the linear algebraic group $\GL_n(F)$ where
 $u^t_n(\centerdot)\in UT_n^t(\centerdot)$,
 $u_n(\centerdot)\in UT_n(\centerdot)$ and
 $d_n(\centerdot)\in D_n(\centerdot)$.
 
 Let $\gl_n(F_R\times  F_L)=
 [u^t_n(F_R)\times u_n(F_L)]\times
 [d _n(F_R)\times d_n(F_L)]$ be the Gauss bilinear decomposition of the matrix $\gl_n(F_R\times  F_L)$ of the algebraic bilinear semigroup
$\GL_n(F_R\times  F_L)$.

Then, if we take into account the maps:
\Bi
\item \quad $u_n(F)\to u_n(F_L)$,
\item \quad $u^t_n(F)\to u^t_n(F_R)$,
\item \quad $d_n(F)\to d_n(F_R)\times d_n(F_L)$,
\Ei
{\bfseries the bilinear algebraic semigroup
{\bbf $\GL_n(F_R\times  F_L)$} is in one-to-one correspondence with the linear algebraic group {\bbf $\GL_n(F)$}\/}.
}
\vskip 11pt

\paragraph{Proof.} \quad
\Be
\item Let $F_R=\{F_{\o\omega _1},\dots,F_{\o\omega _j},
\dots,F_{\o\omega _n}\}$
\resp{$F_L=\{F_{\omega _1},\dots,F_{\omega _j},
\dots,F_{\omega _n}\}$} denote the set of (pseudo-ramified) completions at all archimedean places
$\o\omega _j$ \resp{$\omega _j$} as described in \cite{Pie1}.  Then, to each completion
$F_{\o\omega _j}$ \resp{$F_{\omega _j}$} correspond the conjugacy class representatives
$\gl_n(F_{\o\omega _j})$, $u^t_n(F_{\o\omega _j})$ and $d_n(F_{\o\omega _j})$
\resp{$\gl_n(F_{\omega _j})$, $u_n(F_{\omega _j})$ and $d_n(F_{\omega _j})$} of the respective groups of matrices.

Similarly, to each completion
$F_{\o\omega -\omega _j}=F_{\o\omega _j}\cup F_{\omega _j}$ correspond the conjugacy class representatives
$\gl_n(F_{\o\omega -\omega _j})$,
$u_n(F_{\o\omega -\omega _j})$,
$u^t_n(F_{\o\omega -\omega _j})$ and
$d_n(F_{\o\omega -\omega _j})$ of the respective groups of matrices used in the ``linear'' case.

Then, it is evident that the bilinear algebraic semigroup
$\GL_n(F_R\times F_L)$ is in one-to-one correspondence with the linear algebraic group $\GL_n(F)$ under the three maps of the proposition.

\item Each bilinear conjugacy class representative subspace
$ g^{(n)}(F_{\o\omega -\omega _j})$ of
$\gl_n(F_{\o\omega -\omega _j})$ covers the corresponding linear conjugacy class representative subspace
$g^{(n)}(F_{\o\omega -\omega _j})$ of
$\gl_n(F_{\o\omega -\omega _j})$ if the three maps of the proposition are taken into account.

That is to say there exists the homeomorphism:
\[ H_{g_j^{(n)}}: \quad
g^{(n)}(F_{\o\omega}\times F_{\omega _j})\To
g^{(n)}(F_{\o\omega -\omega _j})\]
in such a way that, if
$g^{(n)}(F_{\o\omega -\omega _j})$ is connected, then 
{\bbf each neighborhood $U_{\o\omega \times \omega _j}$ of a bipoint of
$g^{(n)}(F_{\o\omega_j}\times F_{\omega _j})$ is homeomorphic to the corresponding neighborhood
$H_{g^{(n)}_j}(F_{\o\omega_j}\times F_{\omega _j})
=U_{\o\omega -\omega _j} $ of a point of
$g^{(n)}(F_{\o\omega -\omega _j})$.}

And thus, {\bbf the $n^2$-dimensional representation space $G^{(n)}(F)$ of the linear algebraic group $\GL_n(F)$ coincides with the $n^2$-dimensional representation bisemispace $G^{(n)}(F_R\times F_L)$ of the bilinear algebraic semigroup $\GL_n(F_R\times F_L)$.} \epr
\Ee
\vskip 11pt

\subsection{Bisemischemes}

\subsubsection{Affine bisemischemes}

In analogy with the classical definitions of algebraic geometry \cite{Gro}, \cite{G-R}, \cite{D-G}, \cite{Mum}, let us setup the following definitions introducing affine semischemes.

\Bean
\item {\bbf A \lr semisheaf $\Fs_L$ \resp{$\Fs_R$}}
\vskip 11pt

Let $X_L$ \resp{$X_R$} be a topological semispace restricted to the upper \resp{lower} half space.

A \lr presemisheaf $\Fs_L$ \resp{$\Fs_R$} of abelian semigroups consists of:
\Be
\item[$\alpha $)] for every open \lr subset $U_L\subseteq X_L$ \resp{$U_R\subseteq X_R$}, an abelian semigroup
$\Fs_L(U_L)$ \resp{$\Fs_R(U_R)$}.

\item[$\beta $)] for every inclusion $V_L\subseteq U_L$ \resp{$V_R\subseteq U_R$} of open subsets of $X_L$ \resp{$X_R$}, a morphism of abelian semigroups
\[\res_{U_LV_L}: \quad\Fs(U_L)\to\Fs(V_L)\qquad
\rresp{\res_{U_RV_R}: \quad\Fs(U_R)\to\Fs(V_R)}\;.\]
\Ee

The \lr semisheaf $\Fs_L$ \resp{$\Fs_R$} is a \lr presemisheaf of which sections are determined by local data.

\item A \lr {\bbf ringed semispace} $(X_L,\Fs_L(X_L))$
\resp{$(X_R,\Fs_R(X_R))$} consists of a topological semispace $X_L$ \resp{$X_R$} and of a semisheaf of semirings $\Fs_L$ \resp{$\Fs_R$} on $X_L$ \resp{$X_R$}.

\item  Let $A_L$ \resp{$A_R$} be the polynomial ring $A$ restricted to the ideal
$I_L=\{p_\mu (x_1,\dots,x_n)\}$
\resp{$I_R=\{p_\mu (-x_1,\dots,-x_n)\}$}.

Then, {\bbf the spectrum of $A_L$ \resp{$A_R$}} is the pair
$(\spec A_L,\Fs_L(A_L))$\linebreak \resp{$(\spec A_R,\Fs_R(A_R))$} consisting of the  topological semispace $\spec A_L$
\resp{$\spec A_R$}, being the set of closed subsets associated with all prime ideals of $A_L$ \resp{$A_R$}, together with the semisheaf $\Fs_L( A_L)$ \resp{$\Fs_R( A_R)$} of semirings on it defined as follows:

Let $A_{p_L}$ \resp{$A_{p_R}$} be the localization of $A_L$ \resp{$A_R$} at each prime ideal $p_L\subseteq A_L$
\resp{$p_R\subseteq A_R$}.  For an open set $U_L\subseteq \spec A_L$ \resp{$U_R\subseteq \spec A_R$}, we define the set $\Fs_L(U_L)$ \resp{$\Fs_R(U_R)$} of functions 
$s_L:U_L\to \bigsqcup\limits_{p_L\in A_L} A_{p_L}$
\resp{$s_R:U_R\to \bigsqcup\limits_{p_R\in A_R} A_{p_R}$} such that $s_L(p_L)\in A_{p_L}$ \resp{$s_R(p_R)\in A_{p_R}$} \cite{Hart}.
\Ee

{\bbf An affine \lr semischeme\/} is a locally \lr ringed semispace $(X_L,\Fs_L(X_L))$
\resp{$(X_R,\Fs_R(X_R))$} which is isomorphic to the spectrum $(\spec A_L,\Fs_L(A_L))$ \resp{$(\spec A_R,\Fs_R(A_R))$} of some semiring $A_L$ \resp{$A_R$}.

So, we can introduce {\bbf affine bisemischemes\/} according to the preliminary definitions:
\Bean
\item {\bbf A bisemisheaf $\Fs\RL(X\RL)$ of bilinear semigroups\/} is defined by the tensor product $\Fs_R(X_R) \otimes_{R_R\times R_L}\Fs_L(X_L)$ of the right semisheaf
$\Fs_R(X_R)$ by its left equivalent $\Fs_L(X_L)$ if:
\Be
\item $\Fs\RL$ is defined on the product
$X\RL=X_R\times X_L$ of the right topological semispace $X_R$ by the left topological semispace $X_L$, respectively defined on the right and left semirings $R_R$ and  $R_L$, in such a way that $X_L$ \resp{$X_R$} be a $R_L$-semimodule \resp{$R_R$-semimodule}.

\item $\Fs_L(X_L)$ \resp{$\Fs_R(X_R)$} is a $\Fs_L(\spec A_L)$-semimodule \resp{$\Fs_R(\spec A_R)$-semi\-module}, where $\spec A_L$ \resp{$\spec A_R$} is associated with $X_L$ \resp{$X_R$} \cite{Hart}, in such a way that 
$\Fs_R(X_R) \otimes_{R_R\times R_L}\Fs_L(X_L)$ be a
$\Fs_R(\spec A_R)\times\Fs_L(\spec A_L)$-bisemimodule.
\Ee

Thus, $\Fs_R(X_R) \otimes_{R_R\times R_L}\Fs_L(X_L)$ is a bisemisheaf of bisemimodules over 
$ \spec (A_R) \times \spec (A_L)$ which allows its splitting into:
\[ \Fs_R(X_R)\otimes \Fs_L(X_L)
= (\Fs_R(X_R)\otimes_D \Fs_L(X_L))
\oplus (\Fs_R(X_R)\otimes_{OD} \Fs_L(X_L))\]
according to proposition 2.4.3.

\item {\bbf A ringed bisemispace
$(X_R\times X_L,\Fs_R(X_L)\otimes \Fs_L(X_L))$\/} is then naturally introduced from the product of the right ringed semispace $(X_R,\Fs_R(X_R))$ by its left equivalent\linebreak
$(X_L,\Fs_L(X_L))$.

\item {\bbf The bispectrum of $A_R\times A_L$\/} is introduced by the pair $(\spec A_R\times \spec A_L,\Fs_R(A_R)\otimes \Fs_L(A_L))$ and corresponds to the product of the spectrum $A_R$ by the spectrum of $A_L$.
\Ee

Finally, {\bbf an affine bisemischeme\/} is a locally ringed bisemispace $(X_R\times X_L,\Fs(X_R)\otimes \Fs(X_L))$ which is isomorphic to the spectrum
$(\spec A_R\times \spec A_L,\Fs_R(A_R)\otimes \Fs_L(A_L))$ of the bisemiring $A_R\times A_L$.
\vskip 11pt

\subsubsection[Affine bisemischeme over the bilinear algebraic semigroup $G^{(n)}(F_R\times F_L)$]{\bbf Affine bisemischeme over the bilinear algebraic semigroup $G^{(n)}(F_R\times F_L)$}

\Bean
\item
If the topological \lr semispace is the representation space $G^{(n)}(F_L)\equiv T^{(n)}(F_L)$
\resp{$G^{(n)}(F_R)\equiv T^{(n)}(F_R)$}, over the set of completions $F_L$ \resp{$F_R$}, of the algebraic semigroup of matrices $T_n(F_L)$ \resp{$T^t_n(F_R)$}, then we can introduce {\bbf the \lr semisheaf
$\Fs_L(G^{(n)}(F_L))$
\resp{$\Fs_R(G^{(n)}(F_R))$}\/} over $G^{(n)}(F_L)$ \resp{$G^{(n)}(F_R)$}.

And, {\bbf the bisemisheaf 
$\Fs_R(G^{(n)}(F_R))\otimes\Fs_L(G^{(n)}(F_L))$\/} is defined as the tensor product of the right semisheaf
$\Fs_R(G^{(n)}(F_R))$ by its left equivalent
$\Fs_L(G^{(n)}(F_L))$.

\item {\bbf The affine bisemischeme over the bilinear algebraic semigroup $G^{(n)}(F_R\times F_L)$\/} is a locally ringed bisemispace
$(G^{(n)}(F_R\times F_L),
(\Fs_R(G^{(n)}(F_R))\otimes
\Fs_L(G^{(n)}(F_L)))$ which is isomorphic to the spectrum
$(\spec A_R\times \spec A_L,\Fs_R(A_R)\otimes \Fs_L(A_L))$ of some bisemiring.
\Ee
\vskip 11pt

\subsubsection[Functions on the conjugacy classes of $G^{(n)}(F_L)$ \resp{$G^{(n)}(F_R)$}]{\bbf Functions on the conjugacy classes of $G^{(n)}(F_L)$ \resp{$G^{(n)}(F_R)$}}

Let $G^{(n)}(F_R\times F_L)$ be the $n^2$-dimensional representation bisemispace of the bilinear algebraic semigroup $\GL_n(F_R\times F_L)$.

Let $\{g^{(n)}(F_{\o\omega _{j,m_j}}
\times F_{\omega _{j,m_j}})\}_{j,m_j}$ ($m_j\in\NN$ referring to the multiplicity) be the set of its bilinear conjugacy class representative subspaces in such a way that
 $\{g^{(n)}(F_{\omega _{j,m_j}})\}$
\resp{ $\{g^{(n)}(F_{\o\omega _{j,m_j}})\}$} be the set of \lr linear conjugacy class representative subspaces of
$G^{(n)}(F_L)\equiv T^{(n)}(F_L)$
\resp{$G^{(n)}(F_R)\equiv T^{(n)}(F_R)$}.

Then, we consider the set  $\wh G^{(n)}(F_L)$
\resp{$\wh G^{(n)}(F_R)$} of smooth continuous functions
$\phi ^{(n)}_{G_L}(x_{g_L})$
\resp{$\phi ^{(n)}_{G_R}(x_{g_R})$}, 
$x_{g_L}\in G^{(n)}(F_L)$
\resp{$x_{g_R}\in G^{(n)}(F_R)$}
on $G^{(n)}(F_L)$
\resp{$G^{(n)}(F_R)$} in such a way that
$\wh G^{(n)}(F_L)$
\resp{$\wh G^{(n)}(F_R)$} be partitioned into subsets of functions on the different conjugacy class representatives of
$G^{(n)}(F_L)$
\resp{$G^{(n)}(F_R)$}.
\vskip 11pt

\subsubsection{Proposition}

{\em The set
$\wh G^{(n)}(F_R\times F_L)=
\{ \phi ^{(n)}_{G_{j_R}}(x_{g_{j_R}})
\otimes \phi ^{(n)}_{G_{j_L}}(x_{g_{j_L}})\}_j$ of continuous smooth bifunctions on the bilinear algebraic semigroup
$G^{(n)}(F_R\times F_L)$ is the set of bisections
$\Gamma (\Fs(G^{(n)}(F_R\times F_L)))$ of the bisemisheaf
$ \Fs(G^{(n)}(F_R\times F_L))$ on
$  G^{(n)}(F_R\times F_L)$.
}
\vskip 11pt

\paragraph{Proof.} \quad On every conjugacy class
$g^{(n)}(F_{\o\omega _j}\times F_{\omega _j})$, $1\le j\le r$, of $G^{(n)}(F_R\times F_L)$, there is a set of smooth continuous bifunctions from
$g^{(n)}(F_{\o\omega _j}\times F_{\omega _j})$ into $\CC$ in such a way that 
$\Fs(g^{(n)}(F_{\o\omega _j}\times F_{\omega _j}))$ is a subbisemisheaf.

And, the set of smooth continuous bifunctions on the set of conjugacy classes of
$G^{(n)}(F_R\times F_L)$ forms a bisemisheaf according to section 3.2.1.\epr
\vskip 11pt

\subsubsection{Corollary}

{\em The set
$\wh G^{(n)}(F_R\times F_L)$ of continuous smooth bifunctions on the bilinear algebraic semigroup
$G^{(n)}(F_R\times F_L)$ forms a bisemialgebra in such a way that
$\wh G^{(n)}(F_L)$
\resp{$\wh G^{(n)}(F_R)$} be the \lr semialgebra on 
$G^{(n)}(F_L)$
\resp{$G^{(n)}(F_R)$}.
}
\vskip 11pt

\paragraph{Proof.} \quad This follows from the definition of bisemialgebras in section 2.5.1.\epr
\vskip 11pt

\subsubsection[The bisemialgebra $L\RL^{1-1}(G^{(n)}(F_R\times F_L))$]{\bbf The bisemialgebra $L\RL^{1-1}(G^{(n)}(F_R\times F_L))$}

The set of all continuous measurable \lr functions on
$G^{(n)}(F_L)$
\resp{$G^{(n)}(F_R)$} satisfying:
\[ \int_{G^{(n)}(F_L)}\L|\phi ^{(n)}_{G_L}(x_{g_L})\R|\ dx_{g_L}<\infty \qquad
\rresp{\int_{G^{(n)}(F_R)}\L|\phi ^{(n)}_{G_R}(x_{g_R})\R|\ dx_{g_R}<\infty}\]
with respect to a Haar measure on 
$G^{(n)}(F_L)$
\resp{$G^{(n)}(F_R)$} is the \lr semialgebra
$L^1_L(G^{(n)}(F_L))$
\resp{$L^1_R(G^{(n)}(F_R))$}.

And, the set of all continuous measurable bifunctions on
$G^{(n)}(F_R\times F_L)$ satisfying:
\[ \int_{G^{(n)}(F_R\times F_L)}
\L|
\phi ^{(n)}_{G_R}(x_{g_R})
\otimes \phi ^{(n)}_{G_L}(x_{g_L})
\R|\ dx_{g_R}\ dx_{g_L}<\infty \]
is the bisemialgebra $L\RL^{1-1}(G^{(n)}(F_R\times F_L))$.
\vskip 11pt

\subsubsection{Proposition}

{\em Let
\begin{align*} p_L: \quad G^{(n)}_R(F_R)&\To G^{(n)}_{R(P)}(F_{R(P)}) \\
 \rresp{p_R: \quad G^{(n)}_L(F_L)&\To G^{(n)}_{L(P)}(F_{L(P)})}
 \end{align*}
 be the linear projective morphism mapping
$G_R^{(n)}(F_R)\equiv G^{(n)}(F_R)$
\resp{$G_L^{(n)}(F_L)\equiv G^{(n)}(F_L)$} onto
$G_L^{(n)}(F_L)$
\resp{$G_R^{(n)}(F_R)$} according to proposition 2.4.6 and let
\[ B_L: \quad G^{(n)}_{R(P)}(F_{R(P)})\To G^{(n)}_{L_R}(F_{L_R})
\qquad
\rresp{B_R: \quad G^{(n)}_{L(P)}(F_{L(P)})\To G^{(n)}_{R_L}(F_{R_L)}}\]
be the bijective linear isometric morphism mapping covariant \resp{contravariant} elements of\linebreak 
$G^{(n)}_{L_R}(F_{R(P)})$ \resp{$G^{(n)}_{L(P)}(F_{L(P)})$} into their respective contravariant \resp{covariant} elements according to proposition 2.4.7.

Then, the composition of morphisms:
\begin{align*}
 B_L\circ p_L: \quad L^{1-1}\RL(G^{(n)}(F_R\times F_L)) &\To L^2_{L\times L}(G^{(n)}(F_L\times F_L))\\
\rresp{B_R\circ p_R: \quad L^{1-1}\RL(G^{(n)}(F_R\times F_L)) &\To L^2_{R\times R}(G^{(n)}(F_R\times F_R))}
\end{align*}
transforms the bisemialgebra
$L^{1-1}\RL(G^{(n)}(F_R\times F_L))$ into the bisemialgebra
$L^2_{L\times L}(G^{(n)}(F_L\times F_L))$
\resp{$L^2_{R\times R}(G^{(n)}(F_R\times F_R))$} which is an extended \lr bilinear Hilbert semispace.
}
\vskip 11pt

\paragraph{Proof.}\quad
\Be
\item The composition of morphisms $B_L\circ p_L$ \resp{$B_R\circ p_R$} transforms the continuous bifunctions
$\phi ^{(n)}_{G_R}(x_{g_R})\otimes \phi ^{(n)}_{G_L}(x_{g_L})\in L^{1-1}\RL(G^{(n)}(F_R\times F_L))$ into the bifunctions
$\phi ^{(n)}_{G_L}(x_{g_L})\otimes \phi ^{(n)}_{G_L}(x_{g_L})\in L^{1-1}_{L\times L}(G^{(n)}(F_L\times F_L))$
\resp{$\phi ^{(n)}_{G_R}(x_{g_R})\otimes \phi ^{(n)}_{G_R}(x_{g_R})\in L^{1-1}_{R\times R}(G^{(n)}(F_R\times F_R))$}  satisfying
\[
\int_{G^{(n)}(F_L\times F_L)}
\L| \phi ^{(n)}_{G_L}(x_{g_L})\R|^2\ dx_{g_L}<\infty 
\qquad
\rresp{\int_{G^{(n)}(F_R\times F_R)}
\L| \phi ^{(n)}_{G_R}(x_{g_R})\R|^2\ dx_{g_R}<\infty}
\] according to propositions 2.4.6 and 2.4.7.

$L^2_{L\times L}(G^{(n)}(F_L\times F_L))$
\resp{$L^2_{R\times R}(G^{(n)}(F_R\times F_R))$}
is the (semi)space of square integrable continuous functions on $G^{(n)}(F_L\times F_L)$
\resp{$G^{(n)}(F_R\times F_R)$} restricted to the upper \resp{lower} half space.

\item It is evident from proposition 2.4.9 that
$L^2_{L\times L}(G^{(n)}(F_L\times F_L))$
\resp{$L^2_{R\times R}(G^{(n)}(F_R\times F_R))$}
 is an extended \lr bilinear Hilbert semispace.\epr
 \Ee
 \vskip 11pt

\subsubsection{Bisemimotives}

Let us remark that, in the philosophy of \BSMST, the Chow pure motives, recalled in section 1.2, become Chow pure bisemimotives consisting of pairs
$(X_R\times X_L,\CH^{2n}(X_R\times X_L))$ where:
\Bi
\item $X_R\times X_L$ is the product of an $n$-dimensional smooth projective right semivariety by its corresponding left equivalent;

\item $\CH^{2n}(X_R\times X_L)\simeq \dfrac{Z^i(X_R\times X_L)}{Z^i_{\rm rat}(X_R\times X_L)}$.
\Ei\vskip 11pt

In this new context of bilinearity, the motives thus acquire their letters patent of nobility, as hoped  by A. Grothendieck.
}
\end{document}